\documentclass[11pt]{amsart}
\usepackage{amscd,amsxtra,amssymb,mathrsfs, bbm}
\usepackage{stmaryrd}

\newtheorem{theorem}{Theorem}[section]
\newtheorem{corollary}[theorem]{Corollary}
\newtheorem{lemma}[theorem]{Lemma}
\newtheorem{proposition}[theorem]{Proposition}

\theoremstyle{definition}
\newtheorem{definition}[theorem]{Definition}
\newtheorem{remark}[theorem]{Remark}
\newtheorem{example}[theorem]{Example}

\DeclareMathAlphabet{\mathpzc}{OT1}{pzc}{m}{it}

\DeclareMathOperator{\Dbcoh}{\mathsf{D^{b}_{coh}}}
\DeclareMathOperator{\D*coh}{\mathsf{D^{*}_{coh}}}

\DeclareMathOperator{\Perf}{\mathsf{Perf}}
\DeclareMathOperator{\Qcoh}{\mathsf{Qcoh}}
\DeclareMathOperator{\Band}{\mathsf{Band}}

\DeclareMathOperator{\rk}{rk}

\DeclareMathOperator{\coker}{\mathsf{coker}}

\renewcommand{\ker}{\mathsf{ker}}
\newcommand{\im}{\mathsf{im}}
\renewcommand{\Im}{\mathsf{Im}}
\renewcommand{\dim}{\mathsf{dim}}
\DeclareMathOperator{\tor}{\mathsf{t}}
\DeclareMathOperator{\Coh}{\mathsf{Coh}}
\DeclareMathOperator{\QCoh}{\mathsf{QCoh}}
\DeclareMathOperator{\VB}{\mathsf{VB}}
\DeclareMathOperator{\TF}{\mathsf{TF}}

\DeclareMathOperator{\Pic}{Pic}

\DeclareMathOperator{\CM}{\mathsf{CM}}
\DeclareMathOperator{\Hot}{\mathsf{Hot}}

\DeclareMathOperator{\Ob}{\mathsf{Ob}}

\DeclareMathOperator{\cT}{\mathsf{T}}
\DeclareMathOperator{\cD}{\mathsf{D}}

\DeclareMathOperator{\krdim}{\mathsf{kr.dim}}

\DeclareMathOperator{\gldim}{\mathsf{gl.dim}}

\DeclareMathOperator{\add}{\mathsf{add}}

\DeclareMathOperator{\Hom}{\mathsf{Hom}}

\DeclareMathOperator{\Ext}{\mathsf{Ext}}

\DeclareMathOperator{\Ann}{\mathsf{ann}}
\DeclareMathOperator{\End}{\mathsf{End}}
\DeclareMathOperator{\Mat}{\mathsf{Mat}}

\input xy
\xyoption{all}

\newcommand{\kk}{k}

\renewcommand{\mod}{\mathsf{mod}}

\newcommand{\kA}{\mathcal{A}}
\newcommand{\kB}{\mathcal{B}}
\newcommand{\kE}{\mathcal{E}}
\newcommand{\kF}{\mathcal{F}}
\newcommand{\kG}{\mathcal{G}}
\newcommand{\kH}{\mathcal{H}}
\newcommand{\kI}{\mathcal{I}}
\newcommand{\kJ}{\mathcal{J}}

\newcommand{\kO}{\mathcal{O}}
\newcommand{\kL}{\mathcal{L}}
\newcommand{\kP}{\mathcal{P}}

\newcommand{\kK}{\mathcal{K}}
\newcommand{\kM}{\mathcal{M}}
\newcommand{\kN}{\mathcal{N}}

\newcommand{\kT}{\mathcal{T}}
\newcommand{\kS}{\mathcal{S}}
\newcommand{\kR}{\mathcal{R}}
\newcommand{\kX}{\mathcal{X}}
\newcommand{\kU}{\mathcal{U}}

\newcommand{\sA}{\mathsf{A}}
\newcommand{\sB}{\mathsf{B}}

\newcommand{\sC}{\mathsf{C}}
\newcommand{\sN}{\mathsf{N}}

\newcommand{\lar}{\longrightarrow}

\newcommand{\rA}{A}
\newcommand{\rB}{B}
\newcommand{\rO}{O}
\newcommand{\rR}{R}

\newcommand{\rQ}{Q}

\newcommand{\idm}{\mathfrak{m}}

\newcommand{\mA}{A}
\newcommand{\mB}{B}
\newcommand{\mM}{M}
\newcommand{\mN}{N}

\newcommand{\mT}{T}
\newcommand{\mE}{E}

\newcommand{\mI}{I}

\newcommand{\mF}{F}

\newcommand{\mX}{X}

\newcommand{\mP}{P}

\newcommand{\mQ}{Q}
\newcommand{\mS}{S}
\newcommand{\mY}{Y}
\newcommand{\mU}{U}

\newcommand{\sX}{X}

\newcommand{\sE}{E}

\newcommand{\sU}{U}

\newcommand{\FF}{\mathbb{F}}
\newcommand{\GG}{\mathbb{G}}
\newcommand{\HH}{\mathbb{H}}

\newcommand{\DD}{\mathbb{D}}
\newcommand{\Ss}{\mathbb{S}}
\newcommand{\TT}{\mathbb{T}}
\newcommand{\II}{\mathbb{I}}
\newcommand{\EE}{\mathbb{E}}
\newcommand{\JJ}{\mathbb{J}}
\newcommand{\PP}{\mathbb{P}}
\newcommand{\ZZ}{\mathbb{Z}}

\newcommand{\cA}{\mathcal{A}}
\usepackage{hyperref}

\begin{document}

\title[Tilting on non-commutative curves]{Tilting
on non-commutative rational projective curves}

\author{Igor Burban}
\address{
Mathematisches Institut,
Universit\"at Bonn,
Endenicher Allee 60,
D-53115 Bonn,
Germany
}
\email{burban@math.uni-bonn.de}

\author{Yuriy Drozd}
\address{
 Institute of Mathematics
National Academy of Sciences,
Tereschenkivska str. 3,
01601 Kyiv, Ukraine, and Max-Plank-Institut f\"ur Mathematik,
Vivatsgasse 7, 53111 Bonn, Germany}

\email{drozd@imath.kiev.ua}
\urladdr{www.imath.kiev.ua/$\sim$drozd}

\dedicatory{Dedicated to Helmut Lenzing
on the occasion of his 70th birthday}

\begin{abstract}
In this article we introduce a  new class of non-commutative  projective curves
and show that in certain cases the derived category of coherent sheaves on them has a tilting complex.
 In particular, we prove  that the right bounded derived category of coherent sheaves on a reduced
rational projective curve with only nodes  and cusps as singularities, can be fully
faithfully embedded into the right bounded derived category of the finite dimensional representations
of a certain finite dimensional  algebra of global dimension two.
  As an application of our approach we show that the dimension
of the bounded derived category of coherent sheaves on a rational projective curve with only nodal
or cuspidal singularities is at most two. In the case of the Kodaira cycles of projective lines,
the corresponding  tilted algebras belong to a well-known class  of gentle algebras.
We work out in details the tilting equivalence
in the case of the Weierstrass nodal curve $zy^2 = x^3 + x^2 z$.
\end{abstract}

\maketitle
\thispagestyle{empty}

\tableofcontents

\section{Introduction}

\noindent
By a result of Beilinson \cite{Beilinson} it is known that the derived category of representations of the Kronecker
quiver
$
\vec{Q}
=
\xymatrix{
\bullet  \ar@/^/[r] \ar@/_/[r]  & \bullet
}
$
over a field $\kk$
is equivalent to the derived category of coherent sheaves on the projective line $\PP^1_\kk$.
This article grew up from an attempt to find a similar ``geometric'' interpretation of
other  finite-dimensional algebras.

Let $\sX$ be a singular  reduced projective curve of arbitrary geometric genus, having only
nodes or cusps as singularities. We introduce a certain sheaf of  $\kO_\sX$--orders
(called the Auslander sheaf)  $\kA = \kA_\sX$ and study the category $\Coh(\kA)$
of coherent left modules on the ringed space $(\sX, \kA_\sX)$.  It turns out that the  global dimension
of $\Coh(\kA)$ is equal to \emph{two} and the  original category of coherent sheaves
$\Coh(\sX)$ can be embedded into $\Coh(\kA)$ in two natural but different ways.
Namely, we construct a pair of  fully faithful functors $\FF, \II: \Coh(\sX) \rightarrow \Coh(\kA)$,
where $\FF$ is right exact and $\II$ is left exact such that their images in $\Coh(\kA)$ are closed
 under extensions.

 The functor $\FF$ has a right adjoint functor $\GG: \Coh(\kA) \rightarrow \Coh(\sX)$, which is
 exact.
 We show that $\GG$
 yields an equivalence between the category $\VB(\kA)$ of the locally projective coherent
 $\kA$--modules and the
 category $\TF(\sX)$ of the torsion  free coherent sheaves on $\sX$.
Moreover, we prove  that  $\Coh(\sX)$ is equivalent to the localization of $\Coh(\kA)$ with respect
to a certain bilocalizing subcategory of torsion $\kA$--modules and $\GG$ can be identified with
the canonical localization functor.

It  turns out that the  derived functor $\mathbb{LF}: D^-\bigl(\Coh(\sX)\bigr) \rightarrow
D^-\bigl(\Coh(\kA)\bigr)
$
 is also fully faithful.
 Moreover, if $\sX$ is a \emph{rational} curve of arbitrary arithmetic genus then  the derived category
 $D^b\bigl(\Coh(\kA)\bigr) $ has a \emph{tilting complex}. The corresponding
 tilted algebra $\Gamma_\sX = \End_{D^b(\kA)}(\kH^\bullet)$ is precisely the algebra described in
 \cite[Appendix A.4]{DGVB}. In particular, $\Gamma_\sX$ has global dimension equal to two and
  the derived categories $D^b\bigl(\Coh(\kA)\bigr)$ and $D^b\bigl(\mod-\Gamma_\sX\bigr)$ are equivalent.
  As a corollary we show  that the dimension of the derived category $D^b\bigl(\Coh(\sX)\bigr)$
 is at most two, confirming a conjecture posed by Rouquier
 \cite{Rouquier}. Combining this derived equivalence with the
 embedding $\mathbb{LF}$, we obtain  an exact  fully faithful functor
 $\Perf(\sX) \rightarrow D^b\bigl(\mod-\Gamma_\sX\bigr)$.

If $\sX$ is either a chain or a cycle of projective lines, then the corresponding
tilted algebra $\Gamma_\sX$ belongs to the class of the so-called \emph{gentle} algebras.
They are known to be of \emph{derived-tame} representation type. This gives an alternative
proof of the tameness of $\Perf(\sX)$, obtained for the first time in \cite{Duke}.
If $\sX$ is a Kodaira cycle  of projective lines, then the image
of the triangulated category $\Perf(\sX)$ inside of $D^b\bigl(\mod-\Gamma_\sX\bigr)$
can be characterized  in a very simple way. Let $$\mathbb{L}\nu: D^b\bigl(\mod-\Gamma_\sX\bigr)
\lar    D^b\bigl(\mod-\Gamma_\sX\bigr)$$ be  the  derived Nakayama functor  then the image
of $\Perf(\sX)$ is equivalent to the full subcategory of the complexes $\mP^\bullet$  satisfying
the property  $\mathbb{L}\nu(\mP^\bullet) \cong \mP^\bullet[1]$.

If  $\sE \subseteq \mathbb{P}^2$
is  a singular Weierstra\ss{} cubic curve over an algebraically closed field
$\kk$ given by the equation
$zy^2 = x^3 + x^2z$, then the corresponding algebra
 $\Gamma_\sE$ is the path algebra
of the following quiver with relations:
$$
\xymatrix
{
\bullet \ar@/^/[rr]^{a} \ar@/_/[rr]_{c} & & \bullet \ar@/^/[rr]^{b} \ar@/_/[rr]_{d} & & \bullet
}
\qquad ba = 0,  \quad dc =  0.
$$
We explicitly describe the objects of $D^b\bigl(\mod-\Gamma_\sE)$ corresponding to
the line bundles of degree zero and to the structure sheaves of the regular points
under the embedding $\Perf(\sE) \rightarrow D^b\bigl(\mod-\Gamma_\sE)$.

We hope that the results of this article  will find applications to the homological mirror symmetry
for degenerations of elliptic curves \cite{Kontsevich} and to the theory of integrable systems,
in particular to the study of solutions of Yang-Baxter equations \cite{Polishchuk1, BK4}.

\medskip
\noindent
\emph{Acknowledgement}. Parts of this work were done during the authors stay at the
Mathematical Research Institute in Oberwolfach within the ``Research in Pairs'' programme from
March 8 -- March 21, 2009.
The research of the first-named author was  supported by the DFG
project Bu-1866/1-2, and the research of the second-named author was
supported by the INTAS grant 06-1000017-9093.
The first-named
author would like to thank Bernhard Keller, Daniel Murfet and Catharina Stroppel for
helpful discussions, and  Martin Kalck for pointing
out  numerous misprints in the previous version of this paper.

\section{Auslander sheaf of orders}\label{S:AuslanderSheaf}

In this section we introduce an interesting class of non-commutative ringed spaces supported
on some algebraic curves over an algebraically closed field $\kk$.
Let $\sX$ be a reduced  algebraic  curve over $\kk$ having only \emph{nodes} or \emph{cusps}
as singularities. This means that for any point $x$ from the singular locus $\mathsf{Sing}(X)$
the completion of the local ring
$\widehat{\kO}_{x}$ is isomorphic to
$\kk\llbracket u,v\rrbracket/uv$ (node) or to $\kk\llbracket u,v\rrbracket/(u^2 - v^3)$ (cusp).

Let $\kO = \kO_\sX$ be the structure sheaf of $X$ and $\kK = \kK_\sX$   the sheaf of rational
functions on $\sX$. Consider the ideal sheaf
$\kI$  of the singular locus $\mathsf{Sing}(X)$ (with respect to the reduced
scheme structure) and the sheaf   $\kF = \kI \oplus \kO$.

\begin{definition}
The \emph{Auslander sheaf} of  a reduced curve $\sX$ with  only nodes or
cusps as singularities is the sheaf of $\kO$--algebras
$\kA =  \kA_\sX := {\mathcal End}_\sX(\kF)$.
\end{definition}

\begin{proposition}
In the notations as above, let $\widetilde\sX \stackrel{\pi}\lar \sX$ be the
normalization of $\sX$ and  $\widetilde\kO := \pi_*(\kO_{\widetilde\sX})$. Then we have:
\begin{itemize}
\item The ideal sheaf $\kI$ of the singular locus of $\sX$ is equal to the \emph{conductor
ideal sheaf}  ${\mathcal Ann}_\sX(\widetilde\kO/\kO)$. In particular,
$\kI$ is an ideal sheaf in $\widetilde\kO$.
\item Moreover, we have: $\kI \cong {\mathcal Hom}_\sX(\widetilde\kO, \kO)$,
$\widetilde\kO \cong {\mathcal End}_\sX(\kI)$ and $$
\kA = {\mathcal End}_\sX(\kF) \cong
\left(
\begin{array}{cc}
\widetilde\kO & \kI \\
\widetilde\kO & \kO \\
\end{array}
\right) \subseteq \Mat_2(\kK) =
\left(
\begin{array}{cc}
\kK & \kK \\
\kK & \kK \\
\end{array}
\right).
$$
In other words, $\kA$ is a sheaf of $\kO$--\emph{orders}.
\end{itemize}
\end{proposition}

\begin{proof}
Let $\kJ = {\mathcal Ann}_\sX(\widetilde\kO/\kO)$ be the conductor ideal. Then the sheaf
$\kO/\kJ$ is supported at the set of points of $\sX$,  where $\sX$ is not normal. Since
$\sX$ is a reduced curve, it is exactly the singular locus of $\sX$, hence $\kJ \subseteq \kI$.   Moreover,
$x \in X$ is either a nodal or cuspidal point, then $\bigl(\kO/\kJ\bigr)_x \cong \kk_x$.
This implies that  $\kJ = \kI$.

By  general properties of the conductor ideal it follows that the morphism of $\kO$--modules
$\kI \rightarrow  {\mathcal Hom}_\sX(\widetilde\kO, \kO)$ mapping a local section
$r \in H^0(\sU, \kI)$ to the morphism   $\varphi_r \in $ \linebreak $\Hom_{\sU}\bigl(\widetilde\kO(U), \kO(U)\bigr)$,
given by the rule $\varphi_r(b) = rb$,  is an isomorphism.
In a similar way, the morphism
$\widetilde\kO \lar {\mathcal End}_\sX(\kI)$,  given on the level of local sections
by  the rule  $b \mapsto \psi_b$, where $\psi_b(a) = ba$,  is an isomorphism too.
\end{proof}

\noindent
In what follows we shall use  the following standard results  on  the category of coherent
sheaves over a sheaf of orders on a quasi-projective algebraic variety.

\begin{theorem}\label{T:genpropoford}
Let $\sX$ be a connected quasi-projective algebraic variety of Krull dimension $n$
over an algebraically closed field
$\kk$  and $\kA$
 be a sheaf of orders on $\sX$.
 Let $\QCoh(\kA)$ be  the category of quasi-coherent
left $\kA$-modules and $\Coh(\kA)$ be the category of coherent left
$\kA$-modules. Then we have:
\begin{enumerate}
\item\label{it1:ofgenpropoford}
The category $\Qcoh(\kA)$ is a locally Noetherian Grothendieck category
and $\Coh(\kA)$ is its subcategory of Noetherian objects.
\item\label{it2:ofgenpropoford} For any quasi-coherent $\kA$--modules $\kH'$ and $\kH''$ we
 have the following local-to-global spectral sequence: $H^p\bigl(\sX,
 {\mathcal Ext}^q_\kA(\kH', \kH'')\bigr) \Longrightarrow
 \Ext^{p+q}_\kA(\kH', \kH'')$.
\item\label{it3:ofgenpropoford} If the variety $\sX$ is projective, then $\Coh(\kA)$ is $\Ext$--finite over $\kk$.
\item\label{it4:ofgenpropoford} Assume the ringed space $(\sX, \kA)$
  has ``isolated singularities'', i.e.~ $\kA_x \cong  \Mat_{n \times
    n}(\kO_x)$ for all but finitely many $x \in \sX$. Then we have:
$$
\gldim\bigl(\QCoh(\kA)\bigr)=\gldim\bigl(\Coh(\kA)\bigr) =
\sup_{x \in \sX}\gldim\bigl(\kA_x-\mod\bigr) =
\sup_{x \in \sX}\gldim\bigl(\widehat{\kA}_x-\mod\bigr),
$$
where
$
\widehat{\kA}_x:= \varprojlim \kA_x/ \idm_x^t \kA_x
$
is the radical completion of the ring $\kA_x$
and $\idm_x$ is the maximal ideal of  $\kO_x$.
\item\label{it5:ofgenpropoford}
The canonical functor
$D^*\bigl(\Coh(\kA)\bigr) \rightarrow  \D*coh\bigl(\QCoh(\kA)\bigr)$ is an equivalence of triangulated categories,
where $\D*coh\bigl(\QCoh(\kA)\bigr)$ is the full subcategory of $D^*\bigl(\QCoh(\kA)\bigr)$ consisting
of  complexes with coherent cohomologies and $* \in \left\{\emptyset, b, +, -\right\}$.
\end{enumerate}
\end{theorem}

\noindent
\emph{Comment on the proof}. (\ref{it1:ofgenpropoford})
We refer to \cite[Chapitre IV]{Gabriel} for the definition and general properties of the locally Noetherian
categories. Let $\kH'$ and $\kH''$ be quasi-coherent
$\kA$--modules. Then we have binfunctorial isomorphisms
$$
\Hom_\kA(\kH', \kH'') \cong H^0\bigl(\sX, {\mathcal Hom}_\kA(\kH', \kH'')\bigr) \cong
\Hom_\sX\bigl(\kO, {\mathcal Hom}_\kA(\kH', \kH'')\bigr).
$$
The spectral sequence (\ref{it2:ofgenpropoford}) is just the spectral sequence of the composition
of two left exact functors. Note that ${\mathcal Ext}^q_\kA(\kH', \kH'')$ is a coherent $\kO$--module
for all $q \ge 0$
provided both sheaves $\kH'$ and $\kH''$ are coherent $\kA$--modules.
If $\sX$ is projective, this implies
that $H^p\bigl(\sX, {\mathcal Ext}^q_\kA(\kH', \kH'')\bigr)$ is finite dimensional over $\kk$.
Hence, $\Coh(\kA)$ is $\Ext$--finite in this case, what proves (\ref{it3:ofgenpropoford}).

\vspace{2mm}
\noindent
Let $\rA$ be a left Noetherian ring. By a result of Auslander \cite[Theorem 1]{Auslander}
$$
\gldim(\rA) = \sup \bigl\{\mathsf{pr.dim}(\rA/\mI) \,| \, \mI \subseteq \rA \,\, \mbox{is a left ideal}
\bigr\}.
$$
In particular,  the global dimension of the
category of all $\rA$--modules is equal to the global dimension of the category of Noetherian
$\rA$--modules. Quite analogous observations show that
$\gldim\bigl(\QCoh(\kA)\bigr)=\gldim\bigl(\Coh(\kA)\bigr)$.
If $(\sX, \kA)$ has isolated singularities, then for any coherent  $\kA$--modules
$\kH'$ and $\kH''$ we have:
$$
\krdim\bigl(\mathsf{Supp}\bigl({\mathcal Ext}_\kA^i(\kH', \kH'')\bigr)\bigr) \le \max\bigl\{0, n-i\bigr\}.
$$
This implies that $\Ext^{n+j}_\kA(\kH', \kH'') \cong
H^0\bigl({\mathcal Ext}^{n+j}_\kA(\kH', \kH'')\bigr)$.
Let $\rA$ be an order over a Noetherian local ring $\rO$ and $r$ be the radical of $\rA$. Then we have:
$$
\gldim(\rA) = \mathsf{pr.dim}(\rA/r)  = \sup \bigl\{m \ge 0 \,| \, \Ext^m_\rA(\rA/r, \rA/r) \ne 0
\bigr\},
$$
just as in \cite[\S\,3]{Auslander} or \cite[Chapter IV.C]{Serre}.
In particular, $\gldim(\rA) = \gldim(\widehat\rA)$.
For (\ref{it5:ofgenpropoford}), we refer to \cite[Theorem VI.2.10 and Proposition VI.2.11]{Borel}, \cite[1.7.11]{KashiwaraSchapira} or \cite{SGA}.

\medskip
\noindent
The following  well-known lemma plays a key role in our theory of non-commutative ringed
spaces.

\begin{lemma}\label{L:wellknown}
Let $\rA$ be a commutative Noetherian ring, $\mT = \rA \oplus \mM$ be a finitely
generated $\rA$-module and  $\Gamma = \End_\rA(\mT)$. Then
the $\rA$-module $\mT$ is a projective left $\Gamma$-module and
the canonical ring homomorphism
$\rA \to \End_\Gamma(\mT)$, $a \mapsto \bigl(t \mapsto at\bigr)$ is an isomorphism.
In other words, $\mT$ has the \emph{double centralizer property}.
\end{lemma}

\begin{proof} First note that $\Gamma$ has a matrix presentation
$$
\Gamma =
\left(
\begin{array}{cc}
\rA & \mM^\vee \\
\mM & \End_\rA(\mM)
\end{array}
\right)
$$
and $1_\Gamma = 1_\rA + 1_\mM = e_\rA + e_\mM$. In particular, $\mT \cong \Gamma\cdot e_\rA$ as a left
projective $\Gamma$-module. Hence,
$
\Hom_\Gamma(\Gamma e_\rA, \Gamma e_\rA) \cong e_\rA \Gamma e_\rA  = \rA,
$
and the canonical morphism $\rA \to \End_\Gamma(\mT)$ is an isomorphism.
\end{proof}

\medskip
\noindent
Next, we shall frequently use the following standard  result from  category theory.

\begin{lemma}\label{L:knowcrit}
Let $\FF: \mathsf{C} \lar \mathsf{D}$ be a functor admitting a right adjoint
 functor $\GG: \mathsf{D} \lar \mathsf{C}$. Let  $\phi: \mathbbm{1}_{\mathsf{C}} \to \GG \circ \FF$ be the
 unit of  the adjunction. Then for any objects $X, Y \in \Ob(\mathsf{C})$ the
following diagram is commutative:
$$
\xymatrix
{
\Hom_{\mathsf{C}}(X, Y) \ar[rr]^\FF \ar[rd]_{(\phi_Y)_*}& & \Hom_{\mathsf{D}}\bigl(\FF(X), \FF(Y)\bigr) \\
  & \Hom_{\mathsf{C}}\bigl(X, \GG \FF(Y)\bigr) \ar[ru]_{\cong} &
}
$$
In particular, $\FF$ is fully faithful if and only if  $\phi$ is an isomorphism of functors.
\end{lemma}

\noindent
The main result of this section is the following theorem.

\begin{theorem}\label{T:main}
Let $X$ be a reduced algebraic curve over $\kk$ having only nodal or cuspidal singularities,
$\kF = \kI \oplus \kO$ and $\kA = {\mathcal End}_\sX(\kF)$.
Then the following properties hold.
\begin{enumerate}
\item\label{it1ofmainthm} The functor $\GG = {\mathcal Hom}_\kA(\kF, \, -\,): \Coh(\kA) \rightarrow \Coh(X)$
is exact. Moreover,
$\FF = \kF \otimes_\kO \,-$ is left adjoint to $\GG$ and
 $\HH = {\mathcal Hom}_\kO(\kF^\vee, \, -\,)$ is right adjoint to $\GG$.
\item\label{it2ofmainthm} We have: $\gldim\bigl(\Coh(\kA)\bigr) = 2$.
\item\label{it3ofmainthm} The functors $\GG$ and $\HH$ yield mutually inverse equivalences between the category
$\VB(\kA)$ of locally projective coherent left $\kA$-modules  and  the category
of coherent torsion free sheaves $\TF(X)$. In particular, the
derived functor $\mathbb{R} \GG: D^b\bigl(\Coh(\kA)\bigr) \lar    D^b\bigl(\Coh(\sX)\bigr)$ is essentially surjective.
\item\label{it4ofmainthm} The unit of the adjunction
$\phi: \mathbbm{1}_{\Coh(\sX)} \to \mathbb{G} \circ  \mathbb{F}$ is an isomorphism.
In particular,  the functor $\mathbb{F}$ is fully faithful.
\item\label{it5ofmainthm} Moreover, the unit of the derived adjunction
$\psi: \mathbbm{1}_{D^-(\Coh(\sX))} \to  \mathbb{RG} \circ  \mathbb{LF}$ is also an isomorphism
and the derived functors
$\mathbb{LF}: D^-\bigl(\Coh(\sX)\bigr) \lar  D^-\bigl(\Coh(\kA)\bigr)$  and
$\mathbb{LF}: \Perf(\sX)  \lar  D^b\bigl(\Coh(\kA)\bigr)$ are
 fully faithful.
\end{enumerate}
\end{theorem}

\begin{proof} (\ref{it1ofmainthm}) First note that $\kF$ is endowed with  a natural
 left module structure over the sheaf of algebras
${\mathcal End}_\sX(\kF)$. Moreover, consider the global sections
\begin{equation}\label{E:idempotents}
e_1 =
\left(
\begin{array}{cc}
1 & 0 \\
0 & 0
\end{array}
\right) \quad
\mbox{and} \quad
e_2 =
\left(
\begin{array}{cc}
0 & 0 \\
0 & 1
\end{array}
\right)
\end{equation}
of the sheaf $\kA$. Then we have: $\kF \cong \kA \cdot e_2$, hence $\kF$ is locally projective, viewed
as a left $\kA$--module.  This implies that the functor $\GG$ is exact.

The fact that the functor $\FF$ is left adjoint to $\GG$ is obvious. To see that $\GG$ possesses a right
adjoint functor, consider the coherent sheaf $\kF^\vee \cong  \widetilde\kO \oplus \kO$ and
the sheaf of algebras $\kB =  {\mathcal End}_\sX(\kF^\vee)$. Since  $\kF$ is a torsion
free coherent $\kO$--module,
it is locally Cohen-Macaulay. Moreover,  $\sX$ is a Gorenstein curve, hence
the contravariant functor $\TF(\sX) \rightarrow \TF(\sX)$, mapping a torsion free sheaf
$\kH$ to $\kH^\vee$,  is
an anti-equivalence. This shows that
$$
\kB = \kA^{\mathsf{op}}  = \left(
\begin{array}{cc}
\widetilde\kO & \widetilde\kO \\
\kI & \kO \\
\end{array}
\right)
$$
and $\kB^{\mathsf{op}} = \kA$. In particular, the category of the right $\kB$-modules
is isomorphic to the category of the left $\kA$-modules. Let $\Coh(\kB^{\mathsf{op}})$ be the
category of coherent \emph{right}  $\kB$-modules,  then we have a functor
${\mathcal Hom}_\sX(\kF^\vee,\,-\,):  \Coh(\sX) \rightarrow  \Coh(\kB^{\mathsf{op}})$  possessing
a left adjoint functor $\kF^\vee \otimes_\kB \,-\,: \Coh(\kB^{\mathsf{op}}) \rightarrow  \Coh(\kA)$.
The sheaf $\kF^\vee$ is locally projective as a right
$\kB$-module. Moreover, one easily sees (for instance, as in
Proposition~\ref{L:keyforcomp}) that $\kF^\vee\simeq{\mathcal Hom}_\kA(\kF,\kA)$,
so we have the following commutative diagram of categories of functors:
$$
\xymatrix
{
\Coh(\kB^{\mathsf{op}}) \ar[rr]^{=}  \ar[dr]_{\kF^\vee \otimes_{\kB} \,-\,}
& & \Coh(\kA) \ar[dl]^{{\mathcal Hom}_\kA(\kF,\,-\,)} \\
 & \Coh(\sX) &
}
$$
Hence, the functor $\HH$ is right adjoint to $\GG$.

\vspace{2mm}

\noindent
(\ref{it2ofmainthm}) By part (\ref{it4:ofgenpropoford}) of Theorem \ref{T:genpropoford}, the
global dimension of $\Coh(\kA)$ can be computed locally.
It is clear that $\kA_x \cong  \Mat_2(\kO_x)$ for a smooth point $x \in \sX$. Note that
$\widehat{\kO_x} \cong  \kk\llbracket  u,v\rrbracket/uv$ if the singularity is a node
and $\kk\llbracket  u,v\rrbracket/(u^2 - v^3)$ if it  is a cusp.
 By a result of Bass \cite[Corollary 7.3]{Bass},
   the indecomposable maximal Cohen-Macaulay modules
over the ring $\rO = \kk\llbracket  u,v\rrbracket/uv$ are $\kk\llbracket  u\rrbracket$, $\kk\llbracket  v\rrbracket$ and $\rO$ itself;
whereas for $\rO = \kk\llbracket  u,v\rrbracket/(u^2 - v^3) \cong \kk\llbracket  t^2, t^3\rrbracket$ they are
$\rO$ and $\kk\llbracket  t\rrbracket$.

Hence, in both cases we have: $\CM(\rO) = \add\bigl(\idm\oplus
\rO\bigr)$. By a result of Auslander and Roggenkamp \cite{AusRog},
the algebra $\widehat{\rA} = \End_\rO(\idm \oplus \rO)$  has global dimension \emph{two}, which proves the claim.
See also Remark \ref{R:Auslanderalg}.

\vspace{2mm}
\noindent
(\ref{it3ofmainthm}) Let $x \in X$ be an arbitrary point of the curve $\sX$, $\rO = \kO_x$ and
$\mF = \kF_x$ and $\rB = \End_\rO(\mF)$. Note that $\mF \cong  \mF^\vee$. Then we have a pair of adjoint functors
$$\bar{\HH} = \Hom_\rO(\mF,\,-\,): \; \mod-\rO \rightarrow  \mod-\rB \; \; \mbox{and} \; \;
\bar{\GG} = F \otimes_\rB \,-\,:  \; \mod-\rB \rightarrow   \mod-\rO.$$
By a result of Bass \cite[Corollary 7.3]{Bass},  any maximal Cohen-Macaulay
$\rO$-module is a direct summand of $\mF^n$ for some $n \ge 0$, so
the functors $\bar{\HH}$ and $\bar{\GG}$ induce mutually inverse equivalences between
the category $\CM(\rO)$  of the maximal Cohen-Macaulay modules and the category
$\mathsf{pro}(\rB)$ of the finitely generated  projective right $\rB$-modules.
In particular, the canonical transformations of functors given by the adjunction
$\bar{\zeta}: \mathbbm{1}_{\mathsf{pro}(\rB)} \to \bar{\HH} \circ \bar{\GG}$ and
$\bar{\xi}: \bar{\GG} \circ \bar{\HH} \to \mathbbm{1}_{\CM(\rO)}$ are isomorphisms.

Since the functors $\GG$ and $\HH$ form an adjoint pair, we have natural transformations
of functors
$\zeta: \mathbbm{1}_{\VB(\kB)} \to {\HH} \circ {\GG}$ and
${\xi}: {\GG} \circ {\HH} \to \mathbbm{1}_{\TF(\sX)}$. For any torsion free sheaf
$\kH$ on the curve $\sX$ we have a morphism of $\kO$-modules
$\xi_\kH: {\GG} \circ {\HH}(\kH)  \to \kH$. Moreover, for any point $x \in X$ we
have: $(\xi_\kH)_x = \bar{\xi}_{\kH_x}$, hence $\xi_\kH$ is an isomorphism for
any $\kH$. This implies that $\xi$ is an isomorphism of functors.
In a similar way, $\zeta$ is an isomorphism of functors, too. Hence, the categories
$\TF(\sX)$ and $\VB(\kA)$ are equivalent.

In order to show that $\GG$ is essentially surjective, first note that
any object in $D^b\bigl(\Coh(\kA)\bigr)$ has a finite resolution
by a complex whose terms are locally projective $\kA$-modules.
Since the locally Cohen-Macaulay $\kO$--modules are precisely the torsion free
$\kO$--modules, any object in $D^b\bigl(\Coh(\sX)\bigr)$ is quasi-isomorphic to a bounded complex
whose terms are torsion free coherent $\kO$--modules. Since $\GG$ establishes an equivalence
between $\VB(\kA)$ and $\TF(\sX)$ and is exact,
any object in $D^b\bigl(\Coh(\sX)\bigr)$ has a pre-image
in  $D^b\bigl(\Coh(\kA)\bigr)$.

\vspace{2mm}

\noindent
(\ref{it4ofmainthm})
Let $\phi: \mathbbm{1}_{\Coh(\sX)} \to \mathbb{G} \circ  \mathbb{F}$ be the
unit of the adjunction. Since $\GG$ is exact
and $\FF$ right exact, the composition $\GG \circ \FF$ is right exact, too.
By Lemma \ref{L:wellknown} we know  that the canonical morphism of
sheaves of $\kO$--algebras $\phi_\kO:\kO\to{\mathcal End}_\kA(\kF)$ is an isomorphism.
It implies  that for any locally free coherent $\kO$--module
$\kE$  the canonical morphism $\phi_\kE: \kE \to \GG\FF(\kE)$
is an isomorphism.

 Let $\kN$ be a coherent sheaf
on $\sX$. Since $\sX$ is quasi-projective, we have a presentation
$\kE_1 \stackrel{f}\to \kE_0 \to \kN \to 0$, where $\kE_0$ and $\kE_1$ are locally free.
This gives  a commutative diagram with exact rows
$$
\xymatrix
{
\kE_1 \ar[d]_{\phi_{\kE_1}} \ar[rr]^f &  & \kE_0  \ar[d]^{\phi_{\kE_0}} \ar[rr] & &
\kN \ar[d]^{\phi_{\kN}} \ar[r] & 0 \\
\mathbb{G} \circ \mathbb{F}(\kE_1)  \ar[rr]^{\mathbb{G} \circ \mathbb{F}(f)} &  &\mathbb{G}\circ \mathbb{F}(\kE_0)
\ar[rr] & &  \mathbb{G} \circ \mathbb{F}(\kN)  \ar[r] & 0,
}
$$
where $\phi_{\kE_0}$ and $\phi_{\kE_1}$ are isomorphisms. Hence,
$\phi_{\kN}$ is an isomorphism for any coherent sheaf $\kN$ on the curve $\sX$ and
 $\phi$ is an isomorphism of functors.
The fact that  $\FF$ is fully faithful follows from Lemma \ref{L:knowcrit}.

\vspace{2mm}
\noindent
(\ref{it5ofmainthm})  The derived functors $\mathbb{LF}$ and $\mathbb{RF}$ form again an adjoint pair,
see for example  \cite[Lemma 15.6]{Keller}.
Since $\GG$ is exact, the unit of the adjunction
$\psi: \mathbbm{1}_{D^-(\Coh(\sX)} \to \mathbb{RG}\circ  \mathbb{LF}$
coincides with $\mathbb{L}\phi$. Since $\phi$ is an isomorphism, the natural
transformation $\mathbb{L}\phi$ is an isomorphism, too. Lemma \ref{L:knowcrit} implies  that
the derived functor
$\mathbb{LF}: D^-\bigl(\Coh(\sX)\bigr) \lar  D^-\bigl(\Coh(\kA)\bigr)$  and
its restriction on the category of perfect complexes
$\mathbb{LF}: \Perf(\sX)  \lar  D^b\bigl(\Coh(\kA)\bigr)$ are
 fully faithful.
\end{proof}

\begin{remark}\label{R:Auslanderalg}
Let $\rO$ be either  $\kk\llbracket  x,y\rrbracket/xy$ or  $\kk\llbracket  x,y\rrbracket/(y^2-x^3)$,
$\rR$ be the normalization of $\rO$ and $I = \idm$ be the maximal ideal of $O$, which in this case
is  the conductor ideal. Consider the following $\rO$--order:
$$\rA = \End_\rO(I \oplus \rO) =
\left(
\begin{array}{cc}
\rR & I \\
\rR & \rO
\end{array}
\right),
$$
which is the \emph{Auslander algebra} of the ring $\rO$.
Then $\rA$ is isomorphic to the completion of the path algebra of the
following quiver with relations:
$$
\xymatrix
{
\stackrel{1}\circ \ar@/^/[rr]^{a_{-}}  & &  \stackrel{2}\circ \ar@/^/[ll]^{a_{+}}
 \ar@/_/[rr]_{b_{+}}
 & &
\ar@/_/[ll]_{b_{-}} \stackrel{3}\circ}  \qquad b_{+} a_{-} = 0, \quad  a_{+} b_{-} = 0
$$
if the singularity is a node   and to the completion of the path algebra
$$
\xymatrix
{
\stackrel{1}\circ  \ar@(ul, dl)_{a} \ar@/^/[rr]^{b_{+}}  & & \stackrel{2}\circ \ar@/^/[ll]^{b_{-}}
} \qquad a^2 = b_{-} b_{+}
$$
if the singularity is a cusp. Since $\gldim(\rA) = \mathsf{pr.dim}\bigl(\rA/\mathsf{rad}(\rA)\bigr)$,
 to compute the global dimension of $\rA$ it suffices to compute the projective dimension
of the simple $\rA$--modules.

\begin{itemize}
\item If $\rA$ is nodal, the projective resolutions of the simple $\rA$--modules are
\begin{itemize} \item $
0 \to \mP_1 \oplus \mP_3 \xrightarrow{(a_+ b_+)} \mP_2 \to \mS_2 \to 0,
$
\item $0 \to \mP_3 \xrightarrow{b_{+}} \mP_2 \xrightarrow{a_{-}} \mP_1 \to \mS_1 \to 0,$
\item $0 \to \mP_1 \xrightarrow{a_{+}} \mP_2 \xrightarrow{b_{-}} \mP_3 \to \mS_3 \to 0. $
\end{itemize}
\item If $\rA$ is cuspidal, the projective resolutions of the simple $\rA$--modules are
\begin{itemize}
\item
$
0 \to \mP_1    \stackrel{b_{-}}\lar \mP_2 \to \mS_2 \to 0,
$
\item  $0 \to \mP_1 \xrightarrow{\left(\begin{smallmatrix}b_{-} \\ -a \end{smallmatrix}\right)}   \mP_2 \oplus \mP_1 \xrightarrow{\left(b_+ \, a\right)} \mP_1 \to \mS_1 \to 0$.
\end{itemize}
\end{itemize}
\end{remark}

\noindent
We conclude this section by the following easy  observation.

\begin{proposition}\label{P:importobs}
Let $\cD$ be the full subcategory of the derived category $D^b\bigl(\Coh(\kA)\bigr)$
consisting of the complexes $\kG^\bullet$ such that for any point
$x \in \sX$ the localization $\kG^\bullet_x \in \Ob\bigl(D^b(\kA_x-\mod)\bigr)$
has a finite projective resolution by objects from $\add(\kF_x)$.
Then $\cD$ is triangulated, idempotent complete and $\FF: \Perf(\sX) \lar \cD$
is an equivalence of categories.
\end{proposition}

\begin{proof} First note that the image of $\Perf(\sX)$ under functor $\FF: D^-\bigl(\Coh(\sX)\bigr)
\rightarrow D^-\bigl(\Coh(\kA)\bigr)$
belongs to $\cD$.
Consider the pair of the natural transformations
$\mathbbm{1}_{D^-(\sX)} \stackrel{\phi}\lar \GG \circ \FF$ and
$\FF \circ \GG \stackrel{\psi}\lar
\mathbbm{1}_{D^-(\kA)}$. By Theorem \ref{T:main},  the natural transformation  $\phi$ is an isomorphism. Moreover, $\psi|_{\cD}$ is an
isomorphism, too.  Hence,  $\FF$ and $\GG$ are quasi-inverse equivalences between
the categories $\Perf(\sX)$ and $\cD$.
\end{proof}

\section{Auslander--Reiten translation and  $\tau$-periodic complexes}

We first recall  the following important
result about the existence of the Serre
functor in the derived category $D^b\bigl(\Coh(\kA)\bigr)$, where $\kA$ is the Auslander sheaf
of orders attached to a reduced projective curve  $\sX$ with at most nodal
and cuspidal singularities.

\begin{theorem}\label{T:SerreFunctor}
Let $\sX$ be a reduced projective curve having only nodes or cusps
as singularities, $\kF = \kI \oplus \kO$ and $\kA = {\mathcal End}_\sX(\kF)$.
Consider the  $\kA$--bimodule
$$\omega_\kA :=  {\mathcal Hom}_\sX(\kA, \omega_\sX)
\cong
\left(
\begin{array}{cc}
\kI  & \kI \\
\widetilde\kO & \kO
\end{array}
\right) \otimes \omega_\sX,
$$ where
$\omega_\sX$ is the canonical sheaf  of $\sX$. Then the
endofunctor $\tau: \kG^\bullet \mapsto \omega_\kA \stackrel{\mathbb{L}}\otimes \kG^\bullet$
is the  Auslander--Reiten translation in the derived category $D^b\bigl(\Coh(\kA)\bigr)$. This
 means that for any pair of objects  $\kG^\bullet, \kH^\bullet$ of   $D^b\bigl(\Coh(\kA)\bigr)$ we have bifunctorial isomorphisms
$$
\Hom_{D^b(\kA)}\bigl(\kH^\bullet, \kG^\bullet\bigr) \cong
\DD\bigl(\Ext^1_{D^b(\kA)}\bigl(\kG^\bullet, \tau(\kH^\bullet)\bigr)\bigr),
$$
where $\DD = \Hom_\kk(\,-\,,\kk)$ is the duality over the base field.
\end{theorem}

\begin{proof}
Since $\sX$ is a Gorenstein curve and the sheaf $\kA$ is Cohen-Macaulay as an $\kO$--module, we have
a quasi-isomorphism of the following complexes of $\kA$--bimodules:
$$\mathbb{R}{\mathcal Hom}_\sX(\kA, \omega_\sX)\cong \omega_\kA.$$ Hence, Theorem \ref{T:SerreFunctor}
is a special case of
\cite[Theorem A.4]{NaeghvdBergh} and
\cite[Proposition 6.14]{YekutieliZhang}.
\end{proof}

\begin{corollary}\label{C:transformtheta}
Let $\sX$ be a Kodaira cycle of projective lines, $\kO \xrightarrow{w} \omega_\sX$
be an isomorphism given by a no-where vanishing regular differential 1-form
$w \in H^0(\omega_\sX)$. Then we have
an injective  morphism of $\kA$--bimodules $\theta = \theta_w: \omega_\kA \rightarrow \kA$
 yielding a natural
transformation of exact endofunctors
$$
\theta: \tau_{D^b(\kA)} \lar  \mathbbm{1}_{D^b(\kA)}
$$
of the derived category $D^b\bigl(\Coh(\kA)\bigr)$. In particular, the category $\cD_\theta$ defined
as the full subcategory
of $D^b\bigl(\Coh(\kA)\bigr)$ consisting of all objects $\kG^\bullet$ such that
$\theta_{\kG^\bullet}$ is an isomorphism, is a  triangulated subcategory
of $D^b\bigl(\Coh(\kA)\bigr)$.
\end{corollary}

\begin{proof} Let $\kG_1^\bullet \xrightarrow{f} \kG_2^\bullet \xrightarrow{g}
\kG_3^\bullet \xrightarrow{h} \kG_1^\bullet[1]$ be a distinguished triangle in
$\Dbcoh(\kA)$. Since $\theta$ is a natural transformation of exact functors, we
have a commutative diagram
$$
\xymatrix
{
\tau(\kG_1^\bullet) \ar[r]^{\tau(f)} \ar[d]_{\theta_{\kG_1^\bullet}} &  \tau(\kG_2^\bullet) \ar[r]^{\tau(g)}
\ar[d]^{\theta_{\kG_2^\bullet}} &
\tau(\kG_3^\bullet) \ar[r]^{\tau(h)} \ar[d]^{\theta_{\kG_3^\bullet}} &  \tau(\kG_1^\bullet)[1] \ar[d]^{\theta_{\kG_1^\bullet[1]}}\\
\kG_1^\bullet \ar[r]^f &  \kG_2^\bullet \ar[r]^g &
\kG_3^\bullet \ar[r]^h &  \kG_1^\bullet[1]. \\
}
$$
This implies that if $\kG_1^\bullet$ and $\kG_2^\bullet$ are objects of $\cD_\theta$,  then
$\kG_3^\bullet$ belongs to $\cD_\theta$, too.
\end{proof}

\noindent
The main goal of this section is to establish other descriptions of the category $\cD_\theta$.  To do this, we
consider the local case first.

\begin{lemma}\label{L:keystat}
Let $\rO$ be a nodal singularity,
$\rR$  its normalization and
$\mI = \Ann_{\rO}(\rR/\rO)$ the conductor ideal. Let
$$
\rA = \End_\rO(\mI \oplus \rO) =
\left(
\begin{array}{cc}
\rR & \mI \\
\rR & \rO
\end{array}
\right),
\,\,
\omega_\rA = \Hom_\rO(\rA, \omega_\rO) =
\left(
\begin{array}{cc}
\mI & \mI \\
\rR & \rO
\end{array}
\right), \,\,
\mF = \rA \cdot e_2 =
\left(
\begin{array}{c}
\mI \\
\rO
\end{array}
\right)
$$
be the Auslander algebra of $\rO$, its dualizing module and the indecomposable  projective module
corresponding to the simple $\rA$--module of projective dimension one.
  Let $\mP^\bullet$ be an object
of the derived category $D^b(\rA-\mod)$,
then the following conditions are equivalent:
\begin{enumerate}
\item\label{it1:ofkeylemma} we have an isomorphism
          $\omega_\rA \stackrel{\mathbb{L}}\otimes_\rA \mP^\bullet \cong \mP^\bullet$.
\item\label{it2:ofkeylemma}  the complex $\mP^\bullet$ is quasi-isomorphic to a bounded
complex of modules with entries from $\mathsf{add}(\mF)$.
\item\label{it3:ofkeylemma}  we have an isomorphism
$(\rA/\omega_\rA) \stackrel{\mathbb{L}}\otimes_\rA \mP^\bullet \cong  0$.
\end{enumerate}
\end{lemma}
\begin{proof}
We first consider the case when $\rO$ is complete. Then $\rO = \kk\llbracket x,y\rrbracket/xy$ and
$\rR = \kk\llbracket x\rrbracket \times \kk\llbracket y\rrbracket$.
By
 \cite[Lemma 6.4.1]{vdBergh1}, the functor $\tau = \omega_\rA
 \stackrel{\mathbb{L}}\otimes_\rA -$ is an auto-equivalence of
 $D^b(\rA-\mod)$. Next, the category $D^b(\rA-\mod)$ is Krull--Schmidt
 and
 there are exactly three indecomposable projective $\rA$--modules:
 $$
 \mP_x =
 \left(
\begin{array}{c}
\kk\llbracket x\rrbracket \\
\kk\llbracket x\rrbracket
\end{array}
\right), \,\,
\mP_y =
 \left(
\begin{array}{c}
\kk\llbracket y\rrbracket \\
\kk\llbracket y\rrbracket
\end{array}
\right)\,\,
\mbox{and} \,\,
\mF =
\left(
\begin{array}{c}
\mI \\
\rO
\end{array}
\right).
 $$
 From the exact sequence
 $
 0 \rightarrow \omega_\rA \xrightarrow{\theta} \rA \rightarrow  \rA/\omega_\rA \rightarrow 0,
 $
 for any complex $\mP^\bullet$ from $D^b(\rA-\mod)$ we get a distinguished triangle
 $$
 \omega_\rA \stackrel{\mathbb{L}}\otimes_\rA \mP^\bullet \stackrel{\theta_{\mP^\bullet}}\lar
 \mP^\bullet \lar (\rA/\omega_\rA) \stackrel{\mathbb{L}}\otimes_\rA \mP^\bullet
 \lar
 \omega_\rA \stackrel{\mathbb{L}}\otimes_\rA \mP^\bullet[1].
 $$
 This implies that  $\tau|_{\mathsf{Hot}^b(\mathsf{add}(\mF))} \stackrel{\theta}\lar
 \mathbbm{1}_{\mathsf{Hot}^b(\mathsf{add}(\mF))}
 $ is an isomorphism of functors. On the other hand, we have:
 $$
 \tau(\mP_x) = \mI_x := \left(
\begin{array}{c}
x\kk\llbracket x\rrbracket \\
\kk\llbracket x\rrbracket
\end{array}
\right)\quad
\mbox{and} \quad
 \tau(\mP_y) = \mI_y := \left(
\begin{array}{c}
y\kk\llbracket y\rrbracket \\
\kk\llbracket y\rrbracket
\end{array}
\right).
 $$
 Note that  we have the following short exact sequences:
 $$
 0 \lar \mP_y \stackrel{y}\lar  \mF \lar  \mI_x \lar 0 \quad  \mbox{and} \quad
 0 \lar \mP_x \stackrel{x}\lar  \mF \lar  \mI_y \lar 0.
 $$
 For a complex $\mP^\bullet$ from $D^b(\rA-\mod)$ we define its {defect}  $d(\mP^\bullet)$
 as follows:
 $$
 d(\mP^\bullet) = \sup\left\{n \in \mathbb{Z} \left| H^n\bigl(\rA/\omega_\rA \stackrel{\mathbb{L}}\otimes_\rA
 \mP^\bullet\bigr) \ne 0 \right. \right\}.
 $$
 In particular, $d(\mP^\bullet) = -\infty$ if and only if $\mP^\bullet \in
 \Ob\bigl(\mathsf{Hot}^b\bigl(\mathsf{add}(\mF)\bigr)\bigr)$. By the definition
 of the functor $\tau$ it is clear that
 $d(\mP^\bullet) \ge d\bigl(\tau(\mP^\bullet)\bigr)$ and $d(\mP^\bullet) =
 d\bigl(\tau(\mP^\bullet)\bigr)$
 if and only if $\mP^\bullet \in
 \Ob\bigl(\mathsf{Hot}^b\bigl(\mathsf{add}(\mF)\bigr)\bigr)$. Since $D^b(\rA-\mod)$ is a Krull--Schmidt category,
 this shows the equivalence $(\ref{it1:ofkeylemma})
 \Longleftrightarrow (\ref{it2:ofkeylemma})$. The equivalence $(\ref{it2:ofkeylemma})
 \Longleftrightarrow (\ref{it3:ofkeylemma})$ easily follows from existence of minimal projective resolutions over $\rA$.

Now we consider the general case, when $\rO$ is not necessary complete. The implications
$(\ref{it2:ofkeylemma}) \Longrightarrow (\ref{it3:ofkeylemma}) \Longrightarrow (\ref{it1:ofkeylemma})$
are clear in this case as well. In order to show the implication
$(\ref{it1:ofkeylemma}) \Longrightarrow (\ref{it2:ofkeylemma})$,
it is sufficient to show that a complex $\mP^\bullet \in \Ob\bigl(D^b(\rA-\mod)\bigr)$ is quasi-isomorphic to a bounded
complex of modules with entries from $\mathsf{add}(\mF)$ if and only if
$\widehat{\mP}^\bullet \in \Ob\bigl(D^b(\widehat{\rA}-\mod)\bigr)$ is quasi-isomorphic to a bounded
complex of modules with entries from $\mathsf{add}(\widehat{\mF})$.

 By a result of Bass, see  \cite[Corollary 7.3]{Bass},
any indecomposable torsion free $\rO$--module is isomorphic either to $\rO$ or to a direct summand
of $\rR^m$ for some $m \ge 1$. Since the category of torsion free $\rO$--modules is equivalent to the
category of projective $\rA$--modules, any indecomposable projective module is isomorphic
either to $\mF$ or to a direct summand of $\mP^m$ for some $m \ge 1$, where $\mP = \rA \cdot e_1$.
We assume all entries of the complex $\mP^\bullet$ are projective. Then we have decompositions:
$\mP^n = \mP^n_1 \oplus \mP^n_2$ for all $n \in \ZZ$, where $\mP^n_1 \in \add(\mF)$ and $\mP^n_2 \in
\add(\mP)$. In the set of complexes of projective modules homotopic to $\mP^\bullet$, consider
a representative $\mQ^\bullet$ with the smallest possible number $\sum_{n \in \ZZ}\bigl(\rk_\rO(\mP^n_2)\bigr)$. Let
$n$ be the biggest index for which $\mP^n_2 \ne 0$. Then we have:
$$
\mQ^\bullet =
\left(
\dots \lar \mQ^{n-1}_1 \oplus \mQ^{n-1}_2
\xrightarrow{\left(\begin{smallmatrix}\alpha_{n-1} & \beta_{n-1} \\ \gamma_{n-1} & \delta_{n-1}
\end{smallmatrix}\right)}    \mQ^{n}_1 \oplus \mQ^{n}_2
\xrightarrow{\left(\begin{smallmatrix}\alpha_{n} & \beta_{n}
\end{smallmatrix}\right)} \mQ^{n+1}_1
\lar \dots
\right).
$$
The category of complexes over $\widehat\rA$ is Krull-Schmidt. If $\widehat{\mQ}^\bullet$ is  quasi-isomorphic to a complex, whose entries belong to $\add(\widehat{\mF})$ then the morphism
$\hat{\delta}_{n-1}: \widehat{\mQ}^{n-1}_2 \rightarrow \widehat{\mQ}^{n}_2$ is surjective.
Since the completion functor is faithfully flat, the morphism ${\delta_{n-1}}: {\mQ}^{n-1}_2 \rightarrow {\mQ}^{n}_2$ is surjective, too. Since both modules are projective, ${\delta_{n-1}}$ is the  projection
on a direct summand. This implies that the complex $\mQ^\bullet$ is homotopic to a complex of the form
$$
 \left(
\dots \lar \mQ^{n-1}_1 \oplus \bar{\mQ}^{n-1}_2
\lar    \mQ^{n}_1
\lar   \mQ^{n+1}_1
\lar \dots
\right).
$$
Contradiction. Hence, $\mQ^\bullet$ belongs to $\Hot^b\bigl(\add({\mF})\bigr)$, as wanted.
\end{proof}

\noindent
From now on, let $\sX$ be a Kodaira cycle of projective lines.
We   fix a no-where vanishing regular differential form $w \in H^0(\omega_\sX)$ identifying
$\omega_\kA$ with a sheaf of two-sided $\kA$-ideals.
Hence,  we have a short exact sequence
of $\kA$--bimodules
$
0 \rightarrow \omega_\kA \rightarrow \kA \rightarrow \kA/\omega_\kA \rightarrow 0
$
and for any complex $\kG^\bullet$ from $\Dbcoh(\kA)$ there is  a distinguished triangle
$$
\tau(\kG^\bullet) \stackrel{\theta_{\kG^\bullet}}\lar \kG^\bullet
\lar (\kA/\omega_\kA) \stackrel{\mathbb{L}}\otimes \kG^\bullet  \lar \tau(\kG^\bullet)[1].
$$
The following theorem is the main result of this section.

\begin{theorem}\label{T:tauperiod}
Let $\sX$ be a Kodaira cycle of projective lines, $\kI$ be the ideal sheaf
of the singular locus of $\sX$, $\kF = \kI \oplus \kO$ and
$\kA = {\mathcal End}_\sX(\kF)$ be the Auslander sheaf of $\sX$.
For an object  $\kG^\bullet$ of the derived category $\Dbcoh(\kA)$ the following conditions are equivalent:
\begin{enumerate}
\item\label{it:tauperiod1} we have: $\tau(\kG^\bullet) \cong \kG^\bullet$;
\item\label{it:tauperiod2}  the morphism $\theta_{\kG^\bullet}$ is an isomorphism;
\item\label{it:tauperiod3} we have: $(\kA/\omega_\kA) \stackrel{\mathbb{L}}\otimes \kG^\bullet \cong 0$;
\item\label{it:tauperiod4} $\kG^\bullet$ is an object of the category $\cD_\theta$ introduced in
Corollary \ref{C:transformtheta};
\item\label{it:tauperiod5} $\kG^\bullet$ is an object of the category $\cD$ introduced in Proposition \ref{P:importobs}.
\end{enumerate}
\end{theorem}

\begin{proof}
The equivalences  $(\ref{it:tauperiod2}) \Longleftrightarrow (\ref{it:tauperiod3}) \Longleftrightarrow
(\ref{it:tauperiod4})$  and the implication
$(\ref{it:tauperiod2}) \Longrightarrow (\ref{it:tauperiod1})$ are  obvious.

\vspace{2mm}
\noindent
Note that $(\kA/\omega_\kA) \stackrel{\mathbb{L}}\otimes \kG^\bullet  \cong 0$ in
$D^b\bigl(\Coh(\kA)\bigr)$ if and only if $(\kA/\omega_\kA)_x \stackrel{\mathbb{L}}\otimes \kG^\bullet_x  \cong 0$ in $D^b(\kA_x-\mod)$ for all $x \in \sX$.
Hence, the  equivalence  $(\ref{it:tauperiod3}) \Longleftrightarrow (\ref{it:tauperiod5})$
follows from Lemma \ref{L:keystat}.
The implication $(\ref{it:tauperiod1}) \Longrightarrow (\ref{it:tauperiod3})$ can be shown in a similar way.
\end{proof}

\noindent
Combining Proposition \ref{P:importobs} and Theorem \ref{T:tauperiod}, we obtain the following corollary.

\begin{corollary}\label{C:tauperiod}
Let $\sX$ be a Kodaira cycle of projective lines and $\kA$ be its Auslander sheaf.
Then the image of the functor $\FF: \Perf(\sX) \rightarrow \Dbcoh(\kA)$
is the category of complexes $\kG^\bullet$ such that $\tau(\kG^\bullet) \cong \kG^\bullet$, where
$\tau = \omega_\kA \stackrel{\mathbb{L}}\otimes \,-\,$ is the Auslander-Reiten translate in  $\Dbcoh(\kA)$.
\end{corollary}

\section{Serre quotients and perpendicular categories}

Let $\sX$ be a  reduced curve over $\kk$ having only nodal singularities
and $\kA$ be its Auslander sheaf. The main
goal of this section is to construct  two \emph{different} but natural
embeddings of   $\Coh(\sX)$ into  $\Coh(\kA)$
such that its image will be closed  under extensions.

\vspace{2mm}

\noindent
In Theorem \ref{T:main} it was shown that the functor $\mathbb{F}: \Coh(\sX) \rightarrow
\Coh(\kA)$ is fully faithful. The next proposition characterizes  the image of the functor
 $\mathbb{F}$.

\begin{proposition}
Let
$\kM$ be a coherent left $\kA$--module. Then
 there exists a coherent $\kO$--module $\kN$ such that
$\mathbb{F}(\kN) \cong \kM$ if and only if $\kM$ has a locally projective  presentation
$$\kF \otimes_\kO \kE_1 \lar  \kF \otimes_\kO \kE_0 \lar  \kM \lar  0,$$
  where
$\kE_0$ and $\kE_1$ are locally free coherent $\kO$--modules.
\end{proposition}

\begin{proof} One direction is clear: if $\kN$ is a coherent $\kO$--module then
it has a locally free presentation
$
\kE_1 \rightarrow \kE_0 \rightarrow \kN \rightarrow 0
$
inducing  a  locally projective presentation
$
\kF \otimes_\kO \kE_1 \rightarrow   \kF \otimes_\kO \kE_0 \rightarrow  \kF \otimes_\kO \kN \rightarrow  0.
$
Other way around, let $\kE_0$ and $\kE_1$ be locally free coherent $\kO$-modules
such that
$$
\kF \otimes_\kO \kE_1 \stackrel{f}\lar  \kF \otimes_\kO \kE_0 \lar \kM \to 0
$$
is an exact sequence of coherent left $\kA$-modules. Since the functor
$\mathbb{F}$ is fully faithful, there exists a morphism
of $\kO$--modules $g: \kE_1 \to \kE_0$ such that $f = \mathbb{F}(g)$.
Put $\kN := \coker(g)$. Then we have a commutative diagram
$$
\xymatrix
{
\kF \otimes_\kO \kE_1 \ar[r]^{\mathbb{F}(g)} \ar[d]_\cong & \kF \otimes_\kO \kE_0 \ar[r]
\ar[d]^\cong & \kF \otimes_\kO \kN
 \ar[r] & 0 \\
\kF \otimes_\kO \kE_1 \ar[r]^{f} & \kF \otimes_\kO \kE_0 \ar[r] & \kM
 \ar[r] & 0,
}
$$
implying  that $\kM \cong \kF \otimes_\kO \kN$.
\end{proof}

\begin{proposition}\label{P:extclos}
The  category $\Im(\mathbb{F})$ is closed  under extensions in $\Coh(\kA)$.
\end{proposition}

\begin{proof} Let $\kM'$ and $\kM''$ be two coherent left $\kA$ modules belonging to the image
of $\mathbb{F}$ and
$$
0 \lar  \kM' \lar  \kM \lar  \kM'' \lar  0
$$
be
an exact sequence in $\Coh(\kA)$. Then $\kM$ also belongs to the image of $\mathbb{F}$.
Indeed, by the assumption there exists coherent $\kO_\sX$-modules $\kN'$ and $\kN''$ such that
$\kM' \cong \kF \otimes_\kO \kN'$ and $\kM'' \cong \kF \otimes_\kO \kN''$. Take any
locally free presentation
$
\kE_1' \to \kE_0' \to \kN' \to 0
$
of the coherent sheaf $\kN'$. By Serre's vanishing theorems (see
\cite[Section III.5]{Hartshorne})
 there exists an ample line bundle
$\kL$ and a natural number $n \gg 0$ such that
\begin{itemize}
\item the evaluation morphism
$f'' = \mathrm{ev}: \Hom_\sX(\kL^{\otimes -n}, \kN'') \stackrel{\kk}\otimes \kL^{\otimes -n}
\rightarrow   \kN''$ is an epimorphism;
\item we have the vanishing:
$H^1\bigl({\mathcal Hom}_\kA(\kF, \kF \otimes_\kO \kN') \otimes \kL^{\otimes n}\bigr) = 0$.
\end{itemize}

\noindent
Set $\kE''_0 := \Hom_\sX(\kL^{\otimes -n}, \kN'') \stackrel{\kk}\otimes \kL^{\otimes -n}$ and observe that
the coherent left $\kA$-module $\kF \otimes \kE''_0$ is locally projective.  Hence,
${\mathcal Ext}^i_\kA\bigl(\kF \otimes_\kO \kE''_0, \, -\,\bigr) = 0$ for $i \ge 1$.
Moreover,  for any coherent left $\kA$-module $\kG$ we have an isomorphism of functors
$\Hom_\kA(\kG, \,-\,) \cong  H^0\bigl({\mathcal Hom}_\kA(\kG, \,-\,)\bigr)$ from the category
of coherent $\kA$-modules to the category of finite dimensional vector spaces over $\kk$, where
${\mathcal Hom}_\kA(\kG, \,-\,): \Coh(\kA) \to \Coh(\sX)$. This induces  the  short exact sequence:
$$
0 \lar  H^1\bigl({\mathcal Hom}_\kA(\kG, \kG')\bigr) \lar
\Ext^1_\kA(\kG, \kG')  \lar  H^0\bigl({\mathcal Ext}^1_\kA(\kG, \kG')\bigr)
$$
coming from the standard
local-to-global spectral sequence $H^p\bigl({\mathcal Ext}^q_\kA(\kG, \kG')\bigr) \Longrightarrow $
\linebreak $
\Ext^{p+q}_\kA(\kG, \kG')$.
Since the $\kA$--module $\kF \otimes_\kO \kE''_0$ is locally projective,
 it  implies that
\begin{align*}
\Ext^1_\kA(\kF \otimes_\kO \kE_0'', \kF \otimes_\kO \kN') &\cong
H^1\bigl({\mathcal Hom}_\kA(\kF \otimes_\kO \kE_0'', \kF \otimes_\kO \kN')\bigr)\cong \\
&\cong
H^1\bigl({\mathcal Hom}_\kA(\kF, \kF \otimes_\kO \kN') \otimes_\kO {\kE_0''}^\vee\bigr) = 0.
\end{align*}
Hence,  we can lift the  epimorphism
$1 \otimes f'': \kF \otimes_\kO \kE_0'' \rightarrow  \kM''$ to a morphism
$\bar{f}'': \kF \otimes_\kO \kE_0'' \rightarrow  \kM$, yielding a morphism
$\bar{f} = (\bar{f}', \bar{f}''):
\kF \otimes_\kO (\kE_0' \oplus \kE_0'') \rightarrow
\kM$:
$$
\xymatrix
{
         &  \kF \otimes_\kO \kE_0' \ar[d]_{1 \otimes f'} \ar@{.>}[rd]^{\bar{f}'}
&                 & \kF \otimes_\kO \kE_0'' \ar[d]^{1 \otimes f''} \ar@{.>}[ld]_{\bar{f}''} \\
0 \ar[r] & \kM' \ar[r] & \kM \ar[r] & \kM'' \ar[r] & 0
}
$$
By 5-lemma, $\bar{f}$ is an epimorphism.
 In a similar way, we
construct a presentation $\kE_1'' \rightarrow \kE_0'' \rightarrow \kN'' \rightarrow 0$,  inducing
a presentation
$
\kF \otimes_\kO (\kE_1' \oplus \kE_1'') \rightarrow \kF \otimes_\kO (\kE_0' \oplus \kE_0'') \rightarrow
\kM \rightarrow  0.
$
\end{proof}
Summing up,
Theorem \ref{T:main} and Proposition \ref{P:extclos} imply that
the category of coherent sheaves $\Coh(\sX)$ is equivalent to a full subcategory of $\Coh(\kA)$, which
is closed under taking cokernels and extensions. It turns out that  the category $\Coh(\sX)$ can be embedded
into $\Coh(\kA)$ in a completely different way.

\medskip
Recall that for any singular point $x \in \sX$ the algebra
$\rA = \widehat{\kA}_x$ is isomorphic to the  completion of the path
algebra of the following quiver with relations:
\begin{equation}\label{E:nodalAuslanderalg}
\xymatrix
{
\stackrel{1}\bullet \ar@/^/[rr]^{a_{-}}  & &  \stackrel{2}\circ \ar@/^/[ll]^{a_{+}}
 \ar@/_/[rr]_{b_{+}}
 & &
\ar@/_/[ll]_{b_{-}} \stackrel{3}\bullet}  \qquad b_{+} a_{-} = 0, \quad  a_{+} b_{-} = 0
\end{equation}
Let $\cT$ be the full subcategory of the category of torsion coherent $\kA$--modules,
supported at the singular locus of $\sX$ and corresponding to the simple $\rA$--modules
labeled by bullets. Then $\cT$ is a semi-simple abelian category. Moreover,
$\cT$ is a Serre subcategory of $\Coh(\kA)$ and for any $\kT', \kT'' \in \Ob(\cT)$ we have:
$\Ext^1_\kA(\kT', \kT'') = 0$.

\begin{remark}
Although the category $\cT$ is semi-simple, the second extension group \linebreak
$\Ext^2_\kA(\kT', \kT'')$ is not necessarily zero for a pair of objects $\kT', \kT'' \in \Ob(\cT)$.
Indeed, for the two simple  $\rA$--modules $\mS_1$ and $\mS_3$ from $\cT$ we have:
$
\Ext^2_\rA(\mS_1, \mS_3) = \kk = \Ext^2_\rA(\mS_3, \mS_1).
$
\end{remark}

\begin{definition} Following Geigle and Lenzing
  \cite{GeiglLenzPerpCat},  the \emph{perpendicular category}
  $\cT^\perp$ of the Serre subcategory
$\cT$ is defined as follows:
$$
\cT^\perp = \left\{
\kG \in \Ob\bigl(\Coh(\kA)\bigr)\, \bigl| \, \bigr. \Hom_\kA(\kT, \kG) = 0 = \Ext^1_\kA(\kT, \kG)
\quad \mbox{for all} \quad \kT \in \Ob(\cT)\right\}.
$$
In particular,  $\cT^\perp$ is closed under taking kernels and extensions
inside of $\Coh(\kA)$.
\end{definition}

\begin{proposition}\label{P:Katjashint}
Let $\cT_{\mathrm{sim}}$ be the set of the simple objects of $\cT$.
The perpendicular category $\cT^\perp$ has the following description:
$$
\cT^\perp = \Bigl\{
\kG \in \Ob\bigl(\Coh(\kA)\bigr)\, \bigl| \, \bigr.
{\mathcal Hom}_\kA(\kT, \kG) = 0 = {\mathcal Ext}^1_\kA(\kT, \kG)
\, \,  \mbox{for all} \, \, \kT \in \cT_{\mathrm{sim}}\Bigr\}.
$$
Next%
\footnote{The first-named author would like to thank Catharina Stroppel
for drawing his  attention to this fact.}\!,
the category $\VB(\kA)$ of locally projective left  $\kA$--modules
is a full subcategory of $\cT^\perp$.
\end{proposition}

\begin{proof}
First note that for any objects $\kG \in \Ob(\Coh(\kA))$,
$\kT \in \Ob(\cT)$ and  $i \in \mathbb{Z}$,  the sheaves
${\mathcal Ext}^i_\kA(\kT, \kG)$  are torsion $\kO$--modules. In particular,
we have the  isomorphisms:
$
\Ext^i_\kA(\kT, \kG) \cong  H^0\bigl({\mathcal Ext}^i_\kA(\kT, \kG)\bigr).
$
Hence, for any $i \in \mathbb{Z}$,  the vanishing of $\Ext^i_\kA(\kT, \kG)$ is equivalent to the vanishing of
${\mathcal Ext}^i_\kA(\kT, \kG)$. By induction  on  the length we get:
$$
\kG \in \Ob(\cT^\perp)  \,\, \Longleftrightarrow \, \,
{\mathcal Hom}_\kA(\kT, \kG) = 0 = {\mathcal Ext}^1_\kA(\kT, \kG)
\,\, \mbox{for all} \,\, \kT \in \cT_{\mathrm{sim}}.
$$
This implies the first part of the statement.

In order to show that  $\VB(\kA)$ is a full subcategory of $\cT^\perp$,  it is sufficient
to prove that for any indecomposable projective $\rA$--module $\mP$ and
any simple module $\mS \in \cT_{\mathrm{sim}}$ we have:
$$
\Hom_\rA(\mS, \mP) = 0 = \Ext^1_\rA(\mS, \mP).
$$
This vanishing easily follows from the explicit form of the projective resolution of $\mS$ given
in Remark \ref{R:Auslanderalg}.
\end{proof}

\noindent
For a coherent $\kA$--module $\kH$,  consider its maximal subobject $\tor_{\cT}(\kH)$
belonging to $\cT$
and the canonical short exact sequence
$
0 \rightarrow  \tor_{\cT}(\kH) \rightarrow  \kH \rightarrow  \bar\kH \rightarrow  0.
$
Next, let  $\bar\kH \rightarrow \kI^\bullet_{\bar\kH}$ be an injective resolution of $\bar\kH$.
Choose a distinguished triangle in $D^b\bigl(\Coh(\kA)\bigr)$  determined by the canonical evaluation morphism of complexes of sheaves $\mathsf{ev}_{\bar\kH}$:
\begin{equation}\label{E:univext2}
\bar\kH \stackrel{u}\lar \widetilde\kH \stackrel{v}\lar
\bigoplus\limits_{\kT \in \cT_{\mathrm{sim}}}
\Hom_{D^b(\kA)}\bigl(\kT, \kI^\bullet_{\bar\kH}[1]\bigr) \otimes_\kk  \kT
\xrightarrow{\mathsf{ev}_{\bar\kH}}  \kI^\bullet_{\bar\kH}[1].
\end{equation}
Obviously, this triangle corresponds to a representative of the class of the universal extension sequence
\begin{equation}\label{E:univext1}
0 \lar \bar\kH \stackrel{u}\lar \widetilde\kH \stackrel{v}\lar \bigoplus\limits_{\kT \in \cT_{\mathrm{sim}}}
\Ext^1_\kA(\kT, \bar\kH) \otimes_\kk  \kT \lar 0.
\end{equation}

\noindent
Inspired by the work of Geigle and Lenzing \cite[Section 2]{GeiglLenzPerpCat}, we get the following result.

\begin{theorem}
The correspondence $\kH \mapsto \widetilde\kH$ can be extended to
a functor $\JJ: \Coh(\kA) \rightarrow \cT^\perp$. This functor $\JJ$ is left adjoint
to the inclusion $\II: \cT^\perp \rightarrow \Coh(\kA)$. Moreover, the functor
$\cT^\perp \rightarrow \Coh(\kA)/\cT$ defined as the composition
$
\cT^\perp \stackrel{\II}\lar \Coh(\kA) \stackrel{\PP}\lar \Coh(\kA)/\cT
$
is an equivalence of categories. Here $\Coh(\kA)/\cT$ is the Serre quotient category and $\PP$ is
the corresponding projection  functor. In particular, the perpendicular category
$\cT^\perp$ is abelian and the functor $\JJ$ is right exact.
\end{theorem}

\begin{proof}
Denote $\kE_{\cT}(\bar\kH) := \oplus_{\kT \in \cT_{\mathrm{sim}}}
\Ext^1_\kA(\kT, \bar\kH) \otimes_\kk  \kT$.
We check first that  for any coherent $\kA$--module $\kH$,  the corresponding $\kA$--module
$\widetilde\kH$ given by the universal extension sequence   (\ref{E:univext1})
 belongs to $\cT^\perp$.
Let $\kS$ be an arbitrary element of $\cT_{\mathrm{sim}}$ and
$f: \kS \rightarrow \widetilde\kH$ be a non-zero morphism.
  Since $\Hom_\kA(\kS, \bar\kH) = 0$,
the morphism $g := v_*(f) = vf$ is non-zero too. Since the category  $\cT$ is semi-simple, the morphism
$g: \kS \rightarrow \kE_{\cT}(\bar\kH)$ can be identified with  the inclusion morphism
of a direct summand. By the definition of the evaluation morphism,
$\bigl(\mathsf{ev}_{\bar\kH}\bigr) \circ g \ne 0$. But on the other hand,
$\bigl(\mathsf{ev}_{\bar\kH}\bigr) \circ g = \bigl(\mathsf{ev}_{\bar\kH}\bigr) \circ v \circ f = 0$. Contradiction.

 Since $\Hom_\kA(\kS, \bar\kH) = \Hom_\kA(\kS, \widetilde\kH)  =
 \Ext^1_\kA\bigl(\kS, \kE_{\cT}(\bar\kH)\bigr) = 0$,  the short exact
 sequence (\ref{E:univext1}) induces an exact sequence
$$
0 \lar  \Hom_\kA(\kS,\kE_{\cT}(\bar\kH)) \lar
\Ext^1_\kA(\kS, \bar\kH) \stackrel{u_*}\lar \Ext^1_\kA(\kS, \widetilde\kH) \lar 0.
$$
Since the category $\cT$ is semi-simple, we have: $\dim_\kk\bigl(\Hom_\kA(\kS,\kE_{\cT}(\bar\kH))\bigr)
= \dim_\kk\bigl(\Ext^1_\kA(\kS, \bar\kH)\bigr)$. From the dimension reasons we conclude that
$\Ext^1_\kA(\kS, \widetilde\kH) = 0$. Hence,  $\widetilde\kH$ belongs to $\cT^\perp$ as stated.
The functor $\cT^\perp \rightarrow \Coh(\kA)/\cT$ is fully faithful by \cite[Chapitre III]{Gabriel}.
Moreover, by \cite[Proposition 2.2]{GeiglLenzPerpCat} this functor is an equivalence of categories. In particular, the category $\cT^\perp$ is  abelian.

Now we check that the assignment $\kH \mapsto \widetilde\kH$ can be extended to a functor
$\Coh(\kA) \rightarrow \cT^\perp$. Let $f: \kH' \rightarrow \kH''$ be a morphism in $\Coh(\kA)$.
It is easy to see that $f$ maps $\tor_{\cT}(\kH')$ to $\tor_{\cT}(\kH'')$.
For any object $\kH$ we fix a representative of the cokernel
$\kH \xrightarrow{w} \kH/\tor_{\cT}(\kH)$.  Then  we
obtain the induced map $\bar{f}: \bar\kH' \rightarrow \bar\kH''$ such that the following diagram is
commutative:
$$
\xymatrix
{
0 \ar[r] & \tor_{\cT}(\kH') \ar[d] \ar[r] & \kH' \ar[r]^{w'} \ar[d]_f & \bar{\kH}' \ar[r] \ar[d]^{\bar{f}} & 0 \\
0 \ar[r] & \tor_{\cT}(\kH'')  \ar[r] & \kH'' \ar[r]^{w''} & \bar{\kH}'' \ar[r] & 0.
}
$$
Moreover, the assignment $f \mapsto \bar{f}$ is functorial: $\overline{gf} = \bar{g} \bar{f}$
and $\bar{\mathbbm{1}}_\kH = \mathbbm{1}_{\bar\kH}$.
Next, functoriality of the evaluation morphism and  axioms of triangulated categories imply  there exists
a morphism $\tilde{f}: \widetilde\kH' \rightarrow \widetilde\kH''$ making the following diagram
commutative:
$$
\xymatrix
{
  \bar\kH' \ar[r]^{u'}  \ar[d]_{\bar{f}}  & \widetilde\kH' \ar[r]^{v'}  \ar@{.>}[d]^{\tilde{f}} &
\kE_{\cT}(\bar\kH') \ar[r]^{\mathsf{ev}_{\bar\kH'}} \ar[d]^{\bar{f}_*} &  \bar\kH'[1] \ar[d]^{\bar{f}[1]}\\
  \bar\kH'' \ar[r]^{u''}  &  \widetilde\kH'' \ar[r]^{v''} & \kE_{\cT}(\bar\kH'')
  \ar[r]^{\mathsf{ev}_{\bar\kH'}} & \bar\kH''.
}
$$
Since $\Hom_\cA(\kS, \widetilde{\kH}'') = 0$ for all $\kS \in \cT_{\mathrm{sim}}$, such a  morphism
$\tilde{f}$ is unique. So, we obtain a functor $\JJ: \Coh(\kA) \rightarrow \cT^\perp$.

It is easy to see that we have an isomorphism of functors $\xi: \JJ \circ \II \rightarrow \mathbbm{1}_{\cT^\perp}$. Moreover, there is a natural transformation
$\zeta: \mathbbm{1}_{\Coh(\kA)} \rightarrow \II \circ \JJ$, where for a coherent sheaf
$\kH$ the morphism $\zeta_{\kH}$ is defined to be the composition
$\kH \xrightarrow{w} \bar\kH \xrightarrow{u} \widetilde\kH$. From the short exact sequences
defining the sheaves $\bar\kH$ and $\widetilde\kH$ we conclude that for any
object $\kH''$ from the perpendicular category  $\cT^\perp$ the morphisms $w'_*$ and $u'_*$ are the isomorphisms.
Hence, the morphism $\Hom_\kA\bigl(\JJ(\kH'), \kH''\bigr) \rightarrow
\Hom_\kA(\kH', \kH'')$ given by the composition
$$
\Hom_\kA(\widetilde\kH', \kH'') \stackrel{w'_*}\lar \Hom_\kA(\bar\kH', \kH'')
\stackrel{u'_*}\lar  \Hom_\kA(\kH', \kH'')
$$
is an isomorphism. This shows that $\JJ$ is left  adjoint to the embedding
$\II: \cT^\perp \rightarrow \Coh(\kA)$. Since the category $\cT^\perp$ is abelian, the functor
$\JJ$ is right exact.
 \end{proof}

\begin{remark}\label{R:bilocaliz} By \cite[Proposition 2.2]{GeiglLenzPerpCat} we also  know that
the Serre subcategory $\cT$ is \emph{localizing}. This means that the canonical functor
$\PP: \Coh(\kA) \rightarrow \Coh(\kA)/\cT$ has a right adjoint functor $\widetilde\JJ:
\Coh(\kA)/\cT \rightarrow \Coh(\kA)$.  Theorem \ref{T:SerreQuotStrResult} implies
$\PP$ has also a left adjoint functor, hence $\cT$ is even a \emph{bilocalizing} subcategory.
Moreover, the image of $\widetilde\JJ$ belongs to
the perpendicular category $\cT^\perp$ and there is an isomorphism of functors
$\widetilde\JJ \circ \PP \cong \JJ$.
\end{remark}

\noindent
Note that the exact functor  $\GG = {\mathcal Hom}_\kA(\kF,\,-\,): \Coh(\kA) \rightarrow \Coh(\sX)$
vanishes on the category $\cT$.  Using the universal property of
the Serre quotient category, we obtain an exact functor
$$
\EE: \Coh(\kA)/\cT \lar \Coh(\sX)
$$
such that $\EE \circ \PP = \GG$. The main result of this section is the following theorem.

\begin{theorem}\label{T:SerreQuotStrResult}
The functor $\EE: \Coh(\kA)/\cT \rightarrow \Coh(\sX)$ is an equivalence of abelian categories. Moreover, the functors
$
\GG \circ \II: \cT^\perp \lar \Coh(\sX) \quad \mbox{and} \quad
\JJ \circ \FF: \Coh(\sX) \lar \cT^\perp
$
are mutually quasi-inverse equivalences of abelian categories.
\end{theorem}

\begin{proof}
By Remark \ref{R:bilocaliz}, the second statement implies the first one.

\vspace{2mm}

\noindent
Next,
recall that we have  two adjoint pairs of functors:
$$
\xymatrix
{
\cT^\perp \ar@/^/[rr]^{\II}  & &  \Coh(\kA)  \ar@/^/[ll]^{\JJ}
 \ar@/^/[rr]^{\GG}  & & \Coh(\sX). \ar@/^/[ll]^{\FF}
 }
$$
Let
$\FF \circ \GG \stackrel{\eta}\lar  \mathbbm{1}_{\Coh(\kA)}, \mathbbm{1}_{\Coh(\sX)}
\stackrel{\phi}\lar \GG \circ \FF,
\mathbbm{1}_{\Coh(\kA)} \stackrel{\zeta}\lar \II \circ \JJ$ and $\JJ \circ \II \stackrel{\xi}\lar
\mathbbm{1}_{\cT^\perp}
$
be the morphisms given by the adjunction. Then the functors $\JJ \circ \FF$ and  $\GG \circ \II$ also form
an adjoint pair, whose unit and counit are the natural transformations:
$$
\mu:  \mathbbm{1}_{\Coh(\sX)} \stackrel{\phi}\lar \GG \FF  \xrightarrow{\GG(\zeta)\FF}
 \GG \II \circ \JJ \FF
\quad
\mbox{and}
\quad
\nu: \JJ \FF \circ \GG \II \xrightarrow{\JJ(\eta) \II} \JJ \II \stackrel{\xi}\lar \mathbbm{1}_{\cT^\perp}.
$$
In order to show that $\JJ \FF$ and $\GG \II$ are mutually inverse equivalences of categories,
it is sufficient to show the morphisms $\mu_\kG$ and $\nu_\kH$ are isomorphisms
for arbitrary objects $\kG \in \Ob\bigl(\Coh(\sX)\bigr)$  and $\kH \in \Ob(\cT^\perp)$.
 From the construction of the functors $\FF, \GG, \II$ and $\JJ$ it is clear  that
the morphism of $\kO_x$--modules $\bigl(\mu_\kG\bigr)_x$ and the morphism of $\kA_x$--modules $\bigl(\nu_\kH\bigr)_x$ are isomorphisms  provided $x$ is a smooth point of the curve $\sX$.

Let $x \in \sX$ be a singular point.
Recall that the completion functor is faithfully flat (see for example \cite[Theorem 10.17]{AtiyahMacDonald}), hence a morphism of $\kO_x$ modules
$\mM \xrightarrow{f} \mN$ is an isomorphism if and only if the morphism $\widehat{\mM} \xrightarrow{\hat{f}}
\widehat{\mN}$ is an isomorphism of $\widehat{\kO}_x$--modules. Let
$\rO = \widehat{\kO}_x = \kk\llbracket u,v\rrbracket/uv$, $\rA$ be the Auslander algebra of $\rO$,
$\mF = \widehat{\kF}_x$ and $\bar{\cT}$ be full subcategory of $\rA-\mod$ whose objects correspond
to the torsion sheaves from the category
$\cT$  supported at $x$. Let $\bar\FF = \mF \otimes_\rO \,-\,: \rO-\mod \rightarrow \rA-\mod$,
$\bar\GG = \Hom_\rA(\mF,\-\,): \rA-\mod \rightarrow \rO-\mod$,
$\bar\II: \bar{\cT} \rightarrow \rA-\mod$ be the inclusion functors and $\bar\JJ$ be the right adjoint
to $\bar\II$. Then we have a commutative diagram of categories and functors
\begin{equation}\label{E:globaltolocal}
\xymatrix
{
\cT^\perp \ar@/^/[rr]^{\II}  \ar[d] & &  \Coh(\kA)  \ar@/^/[ll]^{\JJ} \ar[d]
 \ar@/^/[rr]^{\GG}  & & \Coh(\sX). \ar@/^/[ll]^{\FF} \ar[d] \\
\bar{\cT}^\perp \ar@/^/[rr]^{\bar\II}  & &  \rA-\mod  \ar@/^/[ll]^{\bar\JJ}
 \ar@/^/[rr]^{\bar\GG}  & & \rO-\mod. \ar@/^/[ll]^{\bar\FF}
 }
\end{equation}
The vertical arrows correspond to the composition of the localization functor with the functor of the radical completion.
Thus, we have to show the natural transformations of functors
$$
\bar\mu:  \mathbbm{1}_{\rO-\mod} \lar
 \bar\GG \bar\II \circ \bar\JJ \bar\FF
\quad
\mbox{and}
\quad
\bar\nu: \bar\JJ \FF \circ \bar\GG \bar\II \lar
\mathbbm{1}_{\bar{\cT}^\perp}.
$$
are isomorphisms. By Lemma \ref{L:wellknown}, the canonical map
$$\rO = \Hom_\rO(\rO, \rO)  \lar
\Hom_\rA\bigl(\bar\FF(\rO), \bar\FF(\rO)\bigr) =
\Hom_\rA(\mF, \mF)$$ is an isomorphism of algebras.  By Proposition \ref{P:Katjashint}, the module
$\mF$ belongs to  $\bar{\cT}^\perp$. Our next goal is to show $\mF$ is a projective generator
of $\bar{\cT}^\perp$. Indeed, by \cite[Proposition 2.2]{GeiglLenzPerpCat} we know the category
$\bar{\cT}^\perp$ is equivalent to the Serre quotient category $\rA-\mod/\bar{\cT}$.
Let $\bar\PP: \rA-\mod \lar \rA-\mod/\bar{\cT}$ be the canonical functor. Then  $\bar\PP(\mF)$ is a generator of $\rA-\mod/\bar{\cT}$, i.e.~any object in $\rA-\mod/\bar{\cT}$ is a quotient of an object from $\mathsf{add}\bigl(\bar\PP(\mF)\bigr)$.
To prove this, it is sufficient to show that for any projective $\rA$--module $\mQ$ the
object $\bar\PP(\mQ)$ is the quotient of an object from $\mathsf{add}\bigl(\bar\PP(\mF)\bigr)$.
In the notations of Remark \ref{R:Auslanderalg},  we have the following short exact sequence in $\rA-\mod$:
$$0 \lar  \mP_3 \xrightarrow{b_{+}} \mP_2  \xrightarrow{a_{-}} \mP_1 \lar \mS_1 \lar 0,$$
yielding the short exact sequence
 $0 \to \bar\PP(\mP_3) \rightarrow
 \bar\PP(\mP_2) \rightarrow \bar\PP(\mP_1)  \to 0$ in the quotient category
 $\rA-\mod/\cT$. In the same way, we have an exact sequence
 $0 \to \bar\PP(\mP_1) \rightarrow
 \bar\PP(\mP_2) \rightarrow \bar\PP(\mP_3)  \to 0$.

Finally, we check that $\bar\PP(\mF)$ is projective in $\rA-\mod/\bar{\cT}$.
Indeed, assume $\bar\PP(\mX) \stackrel{f}\rightarrow \bar\PP(\mF)$ is an epimorphism in $\rA-\mod/\bar\cT$. By the definition
of the Serre quotient category \cite{Gabriel}, such a morphism is represented by a diagram in $\rA-\mod$
$$
\xymatrix
{
\mY \ar[rr]^g \ar@{_{(}->}[d]_i & & \mQ \\
\mX \ar@{.>}[rr]^f & & \mF \ar@{->>}[u]_p
}
$$
where $i$ is a monomorphism with cokernel belonging to $\bar\cT$, $p$ is an epimorphism whose kernel
belongs to $\bar\cT$ and $g$ is a morphism in $A-\mod$.

Since $\mF$ has no subobjects from $\bar{\cT}$, the morphism $p$ is an isomorphism.
A morphism $\bar\PP(g): \bar\PP(\mY) \rightarrow \bar\PP(\mF)$ is an epimorphism in $\rA-\mod/\bar{\cT}$ if and only if
the cokernel of $g$ belongs  to $\bar\cT$. But $\mF$ has no proper quotients  belonging to
 $\bar\cT$.
Hence, $g$ is an epimorphism in  $\rA-\mod$. Since $\mF$ is projective, the morphism $g$
splits, i.e~there exists $j: \mF \rightarrow \mY$ such that $g j = \mathbbm{1}_\mF$.
But then  $ f \circ \bar\PP(j) = \mathbbm{1}_{\bar{\PP}(\mF)}$, hence $\bar\PP(\mF)$ is projective, as wanted.

This implies that $\mF = \bar\JJ \bar\FF(\rO)$ is a projective generator in $\bar\cT^\perp$, hence
the functor $\bar\JJ \bar\FF: \rO-\mod \rightarrow \bar\cT^\perp$ is a exact.
Moreover, the canonical morphism
$
\Hom_\rO(\rO, \rO) \lar \Hom_\rA(\mF, \mF)
$
is an isomorphism. Hence, $\bar\JJ \bar\FF$ is an equivalence of categories and
its adjoint functor $\bar\GG \bar\II$ is an equivalence, too.
\end{proof}

\begin{remark}
Although  $\cT^\perp$ and $\Coh(\sX)$ are \emph{equivalent} abelian categories
and the functors $\II$ and $\FF$ are fully faithful, the full subcategories
$\II(\cT^\perp)$ and $\FF\bigl(\Coh(\sX)\bigr)$ of the category $\Coh(\kA)$ are \emph{different}.
To show this,  it is sufficient to consider the local situation.
Let $\rO = \kk\llbracket u,v\rrbracket/uv$ and $\rA$ be the corresponding Auslander
algebra. Consider the $\rO$--module $\rO_u = \kk\llbracket u \rrbracket$. It has
a presentation
$
\rO \xrightarrow{v}  \rO \rightarrow \rO_u \rightarrow 0.
$
The functor $\FF$ is right exact, moreover, it induces an equivalence
between the category $\add(\rO)$ and the category $\add(\mP_2) = \add(\mF)$.
This implies that $\mX_u = \FF(\rO_u)$
is given by the presentation
$
\mP_2 \xrightarrow{b_{-} b_{+}} \mP_2 \rightarrow  \mX_u \rightarrow   0.
$
It is then easy to see that $\Hom_\rA(\mS_3, \mX_u) = \kk$. Hence,
$\mX_u$ does not belong to $\bar\cT^\perp$.
\end{remark}

\section{Tilting on rational projective curves with nodal and cuspidal singularities}\label{S:Tilting}

\noindent
Let $\sX$ be a reduced \emph{rational} projective curve with only nodes
or cusps as singularities and $\kA$
be its   Auslander  sheaf of orders
on $\sX$.  The main goal  of this section is to show  that
the derived category $\Dbcoh(\kA)$ has a \emph{tilting complex} and is equivalent to the derived
category of finite dimensional right modules  over a certain  finite dimensional algebra $\Gamma_\sX$.

\subsection{Construction of a tilting complex}
Let $\sX = \displaystyle \bigcup_{i=1}^n \sX_i$, where all components  $\sX_i$ are irreducible
and $\widetilde{\sX} \stackrel{\pi}\lar  \sX$ be the normalization map.
Then $\widetilde\sX = \displaystyle \bigcup_{i=1}^n \widetilde{\sX_i}$, where
$\widetilde{\sX_i} \cong \mathbb{P}^1$ is the normalization
of $\sX_i$ and we have:
$\displaystyle
\widetilde\kO = \bigoplus_{i=1}^n\widetilde{\kO_i},
$
where $\widetilde{\kO_i} = \pi_*(\kO_{\widetilde{\sX_i}})$, $1 \le i \le n$.

\begin{lemma}
In the above notations, consider the locally projective
 $\kA$--module $\kP = \kA \cdot e_1$, where $e_1 \in H^0(\kA)$ is the idempotent
  given by the equation \eqref{E:idempotents}.
Then for any pair of line bundles  $\kL_1$ and $\kL_2$ on the curve $X$  of  the same multidegree
we have: $\kP \otimes_\kO \kL_1 \cong \kP \otimes_\kO \kL_2$.
\end{lemma}

\begin{proof}
Since  $\widetilde{\sX}$ is a union of projective lines,  we have:
$\pi^*\kL_1 \cong \pi^*\kL_2$. The projection formula implies that
$
\pi_* \pi^* \kL_1 \cong \widetilde\kO \otimes_\sX \kL_1 \cong \widetilde\kO \otimes_\sX  \kL_2 \cong
\pi_* \pi^* \kL_2,
$
hence  $\kP \otimes_\kO \kL_1 \cong \kP \otimes_\kO \kL_2$.
\end{proof}

\medskip
\noindent
In what follows, we shall use the notation
$$
\kP = \left(
\begin{array}{c}
\widetilde{\kO}  \\
\widetilde{\kO}
\end{array}
\right)
=
\left(
\begin{array}{c}
\widetilde{\kO}_1  \\
\widetilde{\kO}_1
\end{array}
\right) \oplus
\left(
\begin{array}{c}
\widetilde{\kO}_2  \\
\widetilde{\kO}_2
\end{array}
\right) \oplus \dots \oplus
\left(
\begin{array}{c}
\widetilde{\kO}_n  \\
\widetilde{\kO}_n
\end{array}
\right)
= \kP_1 \oplus \kP_2 \oplus \dots \oplus \kP_n.
$$
For a vector  $\underline{m} = (m_1, m_2, \dots, m_n) \in \ZZ^n$ and a line bundle
$\kL\in \Pic(\sX)$ of multi-degree $\underline{m}$ we denote
$\kP(\underline{m}) = \kP \otimes_\sX \kL \cong
\kP_1(m_1) \oplus \kP_2(m_2) \oplus \dots \oplus \kP_n(m_n).$
For  $\underline{m} = (m, m, \dots, m)$ we shall use the notation: $\kP(\underline{m}) = \kP(m)$.

Let $\kH$ be a coherent left $\kA$--module and $e_1, e_2 \in H^0(\kA)$
 be the  idempotents
 given by (\ref{E:idempotents}).  Then, as $\kO$--module, $\kH$
it splits into the direct sum
$\kH  = e_1 \cdot \kH \oplus e_2 \cdot \kH =   \kH_1 \oplus \kH_2$,
 where $\kH_1$ is an  $e_1 \kA e_1 = \widetilde\kO$--module
 with the induced $\kO$--module structure and
$\kH_2$ is an $e_2 \kA e_2 = \kO$--module. Using these notations,  we write
$\kH \displaystyle =
\left(
\begin{array}{c}
\kH_1 \\
\kH_2
\end{array}
\right).
$
Obviously, a left $\kA$--module $\kH$ is torsion free as an $\kO$--module
if and only if both $\kO$--modules  $\kH_1$ and $\kH_2$ are torsion free.

Next, we shall need the following standard technique from the theory of lattices over orders.
Let $\rO$ be a reduced local ring and  $\rQ = \rQ_1 \times \rQ_2 \times \dots \times \rQ_n$ be its total ring of fractions, where  $\rQ_i$ is a  field for all $1 \le i \le n$. Let $\rA$ be an order over $\rO$, then we have:
\begin{align*}
\rQ(\rA) := \rQ \otimes_\rO \rA &\cong \Mat_{s_1 \times s_1}(\rQ_1)
\times \Mat_{s_2 \times s_2}(\rQ_2) \times \dots \times \Mat_{s_n
  \times s_n}(\rQ_n) \\&:= \rQ_1(\rA) \times \rQ_2(\rA) \times \dots \times \rQ_n(\rA).
\end{align*}
Recall that the ring $\Mat_{s_i \times s_i}(\rQ_i)$ is Morita-equivalent to $\rQ_i$ for all $1 \le i \le n$.
For an $\rA$--module $\mM$ consider the $\rQ(\rA)$--module
$\rQ(\mM) = \rQ \otimes_\rO \mM$.  We say   $\mM$  is   torsion free if the canonical morphism
of $\rA$--modules $\mM \rightarrow \rQ(\mM)$ is injective. In that case, we identify $\mM$ with its image
in $\rQ(\mM)$.

\begin{lemma}\label{L:basicsonorders}
In the notations as above, let $\mM$ and $\mN$ be two Noetherian torsion free $\rA$--modules,
$$
\rQ(\mM) =
\left(
\begin{array}{c}
\rQ_1 \\
\rQ_1 \\
\vdots \\
\rQ_1
\end{array}
\right)^{\oplus m_1} \oplus
\left(
\begin{array}{c}
\rQ_2 \\
\rQ_2 \\
\vdots \\
\rQ_2
\end{array}
\right)^{\oplus m_2}
\oplus \dots \oplus
\left(
\begin{array}{c}
\rQ_n \\
\rQ_n \\
\vdots \\
\rQ_n
\end{array}
\right)^{\oplus m_n}
$$
and
$$
\rQ(\mN) =
\left(
\begin{array}{c}
\rQ_1 \\
\rQ_1 \\
\vdots \\
\rQ_1
\end{array}
\right)^{\oplus l_1} \oplus
\left(
\begin{array}{c}
\rQ_2 \\
\rQ_2 \\
\vdots \\
\rQ_2
\end{array}
\right)^{\oplus l_2}
\oplus \dots \oplus
\left(
\begin{array}{c}
\rQ_n \\
\rQ_n \\
\vdots \\
\rQ_n
\end{array}
\right)^{\oplus l_n}.
$$
Then there is an isomorphism of $\rO$--modules $c: \mS(\mM, \mN) \rightarrow \Hom_\rA(\mM, \mN)$, where
$\mS(\mM, \mN)$ is defined as follows:
$$
\Bigl\{
f= (f_1, f_2, \dots, f_n) \in \Mat_{l_1 \times m_1}(\rQ_1) \times \Mat_{l_2 \times m_2}(\rQ_2) \times
\dots \times \Mat_{l_n \times m_n}(\rQ_n) \Bigl| \Bigr. f\cdot \mM \subseteq \mN
\Bigr\}
$$
and $c(f) \cdot m = f\cdot m$ for $m \in \mM$ and $f \in \mS(\mM, \mN)$.
\end{lemma}

\begin{proof}
First note that the morphism $c$ is well-defined and injective. To prove surjectivity, note that
$\rQ(\rA)$ is injective as $\rA$--module. By \cite[Theorem 1]{Auslander} it is sufficient to show
that for any exact sequence $0 \rightarrow \mI \xrightarrow{\alpha} \rA$ and any morphism $\mI
\xrightarrow{\beta}
\rQ(\rA)$ there exists a morphism $\gamma: \rA \rightarrow \rQ(\rA)$
such that $\gamma \beta = \alpha$. To prove it note that $\rQ
\otimes_\rO \rQ \cong \rQ$ and $\rQ \otimes_\rO \,-\,$ is an exact
functor,  hence the morphism $\alpha$ factorizes through $\rQ(\mI)$
and we have a commutative diagram
$$
\xymatrix
{
 \mI \ar@{^{(}->}[r]^\alpha \ar[d]  \ar@/_20pt/[dd]_{\beta} & \rA \ar[d] \\
 \rQ(\mI) \ar@{^{(}->}[r] \ar[d] & \rQ(\rA) \ar@{.>}[ld]\\
 \rQ(\rA)  &
}
$$
Since the category of $\rQ$--modules is semi-simple, we get the factorization we need.

Hence, for any $\rA$--module $\mN$ the module $\rQ(\mN)$ is injective
over $\rA$.  Let $\mM \xrightarrow{g} \mN$ be a morphism of
$\rA$--modules. By the injectivity of $\rQ(\mN)$ there exists a
morphism $\rQ(\mM) \xrightarrow{\bar{g}} \rQ(\mN)$ making the following
diagram commutative:
$$
\xymatrix
{
\mM \ar[rr] \ar[d]_g & & \rQ(\mM) \ar@{.>}[d]^{\bar{g}}\\
\mN \ar[rr]  & & \rQ(\mN).
}
$$
But then $\bar{g} \in \mS(\mM, \mN)$ and $g = c(\bar{g})$.
\end{proof}

Our next aim is to transfer this technique to the case of sheaves of $\kA$--modules, where
$\kA$ is the Auslander order attached to a projective curve with only nodal or cuspidal singularities.

Let $\kK_i$ be the sheaf of rational functions on the irreducible component
 $\sX_i$ and $\kK$ be the sheaf of rational
functions on $\sX$. Then we have:
$\kK \cong \kK_1 \times \kK_2 \times \dots \times \kK_n$. Let $\rQ_i$ be the field of
rational functions on the component $\sX_i$, then the category of coherent $\kK_i$--modules
is equivalent to the category of finite dimensional vector spaces over $\rQ_i$. Let $\kA$ be the Auslander sheaf  of $\sX$ and $\kH$ be
a torsion free coherent $\kA$--module. Then the canonical morphism of $\kA$--modules
$$
\kH  =
\left(
\begin{array}{c}
\kH' \\
\kH''
\end{array}
\right)
\lar \kK(\kH) :=
\kK \otimes_\kO \kH =
\left(
\begin{array}{c}
\kK_1 \\
\kK_1
\end{array}
\right)^{m_1} \oplus
\left(
\begin{array}{c}
\kK_2 \\
\kK_2
\end{array}
\right)^{m_2} \oplus
\dots \oplus
\left(
\begin{array}{c}
\kK_n \\
\kK_n
\end{array}
\right)^{m_n}
$$
is a monomorphism. In what follows, we consider a torsion free $\kA$--module
$\kH$ as a submodule of $\kK(\kH)$.

\begin{proposition}\label{L:keyforcomp}
In the notations as above, let $\kG$ be a torsion free $\kA$--module and
$$
 \kK(\kG)  =
\left(
\begin{array}{c}
\kK_1 \\
\kK_1
\end{array}
\right)^{l_1} \oplus
\left(
\begin{array}{c}
\kK_2 \\
\kK_2
\end{array}
\right)^{l_2} \oplus
\dots \oplus
\left(
\begin{array}{c}
\kK_n \\
\kK_n
\end{array}
\right)^{l_n}.
$$
Consider the sheaf $\kS(\kH, \kG)$ associated with the following presheaf:
$$
\sU \mapsto \Bigl\{f = \bigl(f_{1},  \dots, f_{n}\bigr) \in
\Mat_{l_1 \times m_1}\bigl(\kK_{1}(\sU)\bigr)  \oplus \dots \oplus \Mat_{l_n \times m_n}\bigl(\kK_{n}(\sU)\bigr)  \, \Big| \, f \cdot \kH(\sU) \subseteq \kG(\sU)
\Bigr\}.
$$
Then the canonical morphism of $\kO$--modules $c: \kS(\kH, \kG) \rightarrow {\mathcal Hom}_\kA(\kH, \kG)$
is an isomorphism.
\end{proposition}

\begin{proof}
First note the morphism $c$ is well-defined. Moreover, for any point $x \in \sX$  its stalk
$c_x$ coincides
with the morphism  from Lemma \ref{L:basicsonorders} applied to the $\kO_x$--order $\kA_x$.
Hence, $c$ is an isomorphism of $\kO$--modules.
\end{proof}

\begin{corollary}\label{C:keytool}
Let $\kG$ and $\kH$ be a pair of coherent torsion free left $\kA$--modules such that
$\kG$ is  locally projective. Then for any $n \ge 0$ we have an isomorphism of vector spaces:
$$
\Ext^n_\kA(\kG, \kH) = H^n\bigl(\kS(\kG, \kH)\bigr).
$$
\end{corollary}

\begin{proof}
Indeed, since ${\mathcal Ext}_\kA^i(\kG, \kH) = 0$ for all $i \ge 1$, the local-to-global
spectral sequence implies that $\Ext^i_\kA(\kG, \kH) \cong H^i\bigl({\mathcal Hom}_\kA(\kG, \kH) \bigr)
\cong H^i\bigl(\kS(\kG, \kH) \bigr).$
\end{proof}

\begin{corollary}
Let $\sX$ be a curve with nodal or cuspidal singularities, $\kA$ be its Auslander sheaf, $\kG$ and
$\kH$ be two torsion free coherent $\kA$--modules such that
$$
\kK(\kG)  =
\left(
\begin{array}{c}
\kK_1 \\
\kK_1
\end{array}
\right) \oplus
\left(
\begin{array}{c}
\kK_2 \\
\kK_2
\end{array}
\right) \oplus
\dots \oplus
\left(
\begin{array}{c}
\kK_n \\
\kK_n
\end{array}
\right) \cong \kK(\kH).
$$
Then the $\kO$--module ${\mathcal Hom}_\kA(\kG, \kH)$ is isomorphic to the sheaf
associated with the presheaf
$$
\sU \mapsto \Bigl\{f = \bigl(f_{1},  \dots, f_{n}\bigr) \in
\kK_{1}(\sU)  \oplus \dots \oplus \kK_{n}(\sU)  \, \Big| \,
\begin{array}{l}
 f \cdot  \kH'(\sU) \subseteq \kG'(\sU) \\
 f \cdot \kH''(\sU) \subseteq \kG''(\sU)
\end{array}
\Bigr\}.
$$
In particular, for the locally projective $\kA$--modules $\kP = \kA \cdot e_1$ and $\kF = \kA \cdot
e_2$, where $e_1, e_2 \in H^0(\kA)$ are the idempotents given by \eqref{E:idempotents}, we have:
$$
\widetilde\kO \cong {\mathcal Hom}_\kA(\kP, \kP), \,\,
\kO \cong {\mathcal Hom}_\kA(\kF, \kF), \,\,
\widetilde\kO \cong {\mathcal Hom}_\kA(\kF, \kP)\,\, \mbox{and} \,\,
\kI \cong {\mathcal Hom}_\kA(\kP, \kF).
$$
\end{corollary}

\medskip
\noindent
The following  proposition  plays the key role in our approach to non-commutative
rational projective curves.

\begin{proposition}\label{T:crucial}
Let $\sX$ be a rational reduced projective curve with only nodal or cuspidal singularities
and $\kA$ be its Auslander sheaf of orders.  Consider the torsion $\kA$--module $\kS$
given by its locally free resolution
$$
0 \lar \left(
\begin{array}{c}
\kI \\
\kI
\end{array}
\right)
\lar
\left(
\begin{array}{c}
\kI \\
\kO
\end{array}
\right)
\lar
\kS \lar 0.
$$
Note that the first term of this short exact sequence is isomorphic to $\kP(-2)$, the middle term
is $\kF$ and $\kS = \kS_1 \oplus \kS_2 \oplus \dots \oplus \kS_t$, where the torsion
module $\kS_i$ is supported at the singular point $x_i \in \sX$ and corresponds to the unique
simple $\widehat{\kA}_{x_i}$--module of projective dimension one.
Then the complex
\begin{align*}
\kH^\bullet &:= \kS[-1] \oplus \kP(-1) \oplus \kP =\\ &\ = \bigl(\kS_1 \oplus \kS_2 \oplus \dots \oplus \kS_t\bigr)[-1]
\oplus \bigl(\kP_1(-1) \oplus \kP_1\bigr) \oplus \bigl(\kP_2(-1) \oplus \kP_2\bigr) \oplus
\dots \bigl(\kP_n(-1) \oplus \kP_n\bigr)
\end{align*}
is \emph{rigid}
 in the
derived category of coherent sheaves $D^b\bigl(\Coh(\kA)\bigr)$, i.e. for all $i \ne 0$ we have
$$
\Hom_{D^b(\kA)}\bigl(\kH^\bullet, \kH^\bullet[i]\bigr) = 0.
$$
\end{proposition}

\begin{proof}
By Corollary \ref{C:keytool} we have:
\begin{align*}
&\Ext^i_\kA\bigl(\kP(-1) \oplus \kP,\, \kP(-1) \oplus \kP\bigr) =
H^i\Bigl(\kS\bigl(\kP(-1) \oplus \kP, \kP(-1) \oplus \kP\bigr)\Bigr) =\\
\ &= H^i\bigl(\sX, \tilde\kO(-1) \oplus \tilde\kO^{\oplus 2} \oplus \tilde\kO(1)\bigr) =
H^i\bigl(\widetilde\sX, \kO_{\widetilde\sX}(-1) \oplus
\kO_{\widetilde\sX}^{\oplus 2} \oplus \kO_{\widetilde\sX}(1)\bigr) = 0
\end{align*}
for all $i \ne 0$.

Now  we check the torsion $\kA$--module $\kS$ is also exceptional. Again, using
the local-to-global spectral sequence, we have:
$
\Ext^i_\kA(\kS, \kS) =  H^0\bigl(\sX, {\mathcal Ext}^i_\kA(\kS, \kS)\bigr), \, i \ge 0.
$
Hence, the vanishing of $\Ext^i_\kA(\kS, \kS)$ can be checked locally.
Using the projective resolution of the simple $\widehat{\kA}_{x_i}$--module
$\widehat{\kS}_{x_i}$ given in Remark \ref{R:Auslanderalg}, we get the desired vanishing.

Next,  $\kS$ is torsion and $\kP(n)$ is torsion free, hence we get:
$\Hom_\kA\bigl(\kS, \kP(n)\bigr) = 0$ for all $n \in \mathbb{Z}$. Since
$\kS$ has a locally projective resolution of length one, we have:
${\mathcal Ext}_\kA^i\bigl(\kS, \kP(n)\bigr) = 0$ for $i \ne 1$. The local-to-global
spectral sequence implies that
$
\Ext^i_\kA\bigl(\kS, \kP(n)\bigr) = 0
$
for all $i \in \ZZ$.
Finally, it remains to note that ${\mathcal Ext}^i_\kA\bigl(\kP(n), \kS\bigr) = 0$
for all $n \in \mathbb{Z}$ and $i \ge 0$, so the local-to-global
spectral sequence implies again
$
\Ext^i_\kA\bigl(\kP(n), \kS\bigr)  = 0.
$
\end{proof}

Let $\cD$ be a triangulated category admitting all set-indexed coproducts. Recall that
an object $\mX \in \Ob(\cD)$ is called \emph{compact} if for an arbitrary family
 $\bigl\{\mY_i\bigr\}_{i \in I}$ of objects of $\cD$  the canonical map
 $$
 \oplus_{i \in I} \Hom_{\cD}(\mX, \mY_i) \lar
 \Hom_{\cD}(\mX, \oplus_{i \in I} \mY_i)
 $$
 is an isomorphism.
 An object $\mX$ compactly generates $\cD$ if it is compact and
 $$
 \mX^\perp := \left\{\mY \in \Ob(\cD) \, \Big| \, \Hom_{\cD}(\mX,
 \mY[n]) = 0 \;  \forall \, n \in \ZZ \right\} = 0.
 $$

\medskip
\noindent
Recall the following result of Keller \cite{DerivingDG}.

\begin{theorem}\label{T:Keller}
Let $\cD$ be an algebraic triangulated category admitting all set-indexed coproducts
and $\mX$ be a compact generator of $\cD$ such that
$\Hom_{\cD}\bigl(\mX, \mX[n]\bigr) = 0$ for all $ n \in \ZZ\setminus \{0\}$.
Let $\Gamma =
\End_{\cD}(\mX)$ and $\mathsf{Mod}-\Gamma $ be the category of \emph{all} right $\Gamma$--modules.
Then there exists an exact equivalence of triangulated categories
$
\TT: \cD \lar D\bigl(\mathsf{Mod}-\Gamma\bigr)
$
such that for an arbitrary object $\mY \in \Ob(\cD)$ we have:
$H^n\bigl(\TT(\mY)\bigr) = \Hom_{\cD}(\mX, \mY[n])$, where $\Hom_{\cD}(\mX, \mY[n])$ is endowed
with the natural structure of a right $\Gamma = \End_{\cD}(\mX)$--module.
Such an object $\mX$ is called \emph{tilting} and its endomorphism algebra $\Gamma$ is  the corresponding
\emph{tilted algebra}.
\end{theorem}

\noindent
In order to restrict the equivalence $\TT$ on the derived category of Noetherian objects
of a Grothendieck abelian category, we use the following result of Krause \cite[Proposition 2.3]{Krause}.

\begin{theorem}\label{NeemanKrause}
Let $\sA$ be a locally Noetherian Grothendieck category of \emph{finite global dimension}
and $\mathsf{N}$ be its full subcategory of Noetherian objects. Let $\cD_c(\sA)$ be the category of compact objects
of $D(\sA)$. Then the image of the canonical functor
$D^b(\mathsf{N}) \rightarrow D(\sA)$ is equivalent to $\cD_c(\sA)$.
\end{theorem}

\medskip
\noindent
The following result was explained to the first-named author by Daniel Murfet.

\begin{proposition}\label{P:Murfet}
Let $\sA$ be a locally Noetherian Grothendieck category of finite global dimension, $\sN$ be its full
subcategory of Noetherian objects  and $\cD = D(\sA)$ be the derived category of $\sA$. Let
$\II: D^-(\sN) \rightarrow \cD$ be the canonical functor. Then an object $\mX$ of the category
$\cD$ belongs to the image of $\II$ if and only if for every family of objects
$\{\mY\}_{i \in I}$ of $\cD$ such that $\oplus_{i \in I} \mY_i$ has bounded below
cohomology, the
canonical map $$
 \bigoplus\limits_{i \in I} \Hom_{\cD}\bigl(\mX, \mY_i\bigr) \lar
 \Hom_{\cD}\bigl(\mX, \bigoplus\limits_{i \in I} \mY_i\bigr)
 $$
 is an isomorphism.
\end{proposition}

\begin{proof} We follow the main steps of the proof of \cite[Lemma 4.1]{Krause}.
First check that any object $\mX$ of the category $D^-(\sN)$ has the stated property.
Let $\{\mY_i\}_{i\in I}$ be a family of complexes from $\sA$ with common lower bound for
the non-vanishing cohomology.
Then there exists $n \in \ZZ$ such that for all $j < n$ and $i \in I$ we have:
$H^j(\mY_i) = 0$. Let $\tau_{\ge n}(\mX)$ be the truncation of $\mX$. By the assumption, the complex  $\tau_{\ge n}(\mX)$
has bounded Noetherian cohomology. Since $\sA$ has finite global dimension, the complex
$\tau_{\ge n}(\mX)$ is compact in $\cD$ and we have canonical isomorphisms
$$
\bigoplus\limits_{i \in I} \Hom_{\cD}\bigl(\mX, \mY_i\bigr) \cong \bigoplus\limits_{i \in I}
\Hom_{\cD}\bigl(\tau_{\ge n}(\mX), \mY_i\bigr) \cong
\Hom_{\cD}\bigl(\tau_{\ge n}(\mX), \bigoplus\limits_{i \in I} \mY_i\bigr) \cong
\Hom_{\cD}(\mX, \bigoplus\limits_{i \in I} \mY_i).
$$
Now, let $\mX$ be an object of $\cD$  such that the map $\oplus_{i \in I} \Hom_{\cD}(\mX, \mY_i) \lar
 \Hom_{\cD}(\mX, \oplus_{i \in I} \mY_i)$ is an isomorphism for an arbitrary family of objects
 with a common lower bound for the non-vanishing cohomology. First we check there exists
 $n \in \ZZ$ such that for all $m \ge n$ we have: $H^m(\mX) = 0$.  Let
 $\mX = \bigl(\dots \rightarrow \mX^{j-1} \xrightarrow{\delta^{j-1}} \mX^j
 \xrightarrow{\delta^{j}} \mX^{j+1} \rightarrow \dots\bigr)
$ and assume $H^j(\mX) \ne 0$. Let $\ker(\delta^j) := \bigl(Z^j(\mX), \alpha_j\bigr)$,
$Z^j(\mX) \xrightarrow{\beta_j} H^j(\mX)$ be the canonical epimorphism   and
$H^j(\mX) \xrightarrow{\gamma_j} \mE\bigl(H^j(\mX)\bigr)$
be the injective envelope of $H^j(\mX)$. Note that the composition $\gamma_j \beta_j$ is non-zero.
Moreover, since $\mE\bigl(H^j(\mX)\bigr)$ is injective, there exists a morphism
$\varphi_j: \mX^j \rightarrow \mE\bigl(H^j(\mX)\bigr)$ such that $\varphi_j \alpha_j =
\gamma_j \beta_j$. Next, by the universal property of the kernel, there exists a morphism
$\tilde{\delta}^{j-1}: \mX^{j-1} \rightarrow Z^j(\mX)$ such that $\alpha_j \tilde{\delta}^{j-1} =
\delta^{j-1}$.
$$
\xymatrix
{
\dots \ar[r] & \mX^{j-1} \ar[r]^{\delta^{j-1}} \ar[rd]_{\tilde{\delta}^{j-1}} & \mX^j \ar[r]^{\delta^{j}} \ar@/^3pc/[ddd]^{\phi_j} & \mX^{j+1} \ar[r] & \dots \\
&                                              & Z^j(\mX) \ar@{^{(}->}[u]_{\alpha_j}   \ar@{->>}[d]^{\beta_j}             &                    &       \\
&                                              & H^j(\mX)  \ar@{_{(}->}[d]^{\gamma_j}               &                    &       \\
&                                              & \mE\bigl(H^j(\mX)\bigr)                 &       &
}
$$
Note that $\varphi_j \delta^{j-1} = \gamma_j \beta_j
\tilde{\delta}^{j-1} = 0$. As a result, we get  a morphism $\mX
\rightarrow \mE\bigl(H^j(\mX)\bigr)[-j]$ inducing a non-zero
map in cohomology. Hence, if $\mX$ has unbounded cohomology to the right, the morphism
$\mX \rightarrow \oplus_{j \in \ZZ_+} \mE\bigl(H^j(\mX)\bigr)[-j]$ can not factor
through a finite set of indices in $\ZZ_+$. This implies that the canonical map
$$
 \bigoplus\limits_{j \in \ZZ_+} \Hom_{\cD}\Bigl(\mX, \mE\bigl(H^j(\mX)\bigr)[-j]\Bigr) \lar
 \Hom_{\cD}\Bigl(\mX, \bigoplus\limits_{j \in \ZZ_+} \mE\bigl(H^j(\mX)\bigr)[-j]\Bigr)
 $$
 is not an isomorphism. Contradiction.

 In remains to show that $\mX$ has coherent cohomology. Let $\left\{\mE_i\right\}_{i \in I}$ be an
 arbitrary family of injective objects in $\sA$. Then we have:
 $$
 \Hom_{\cD}\bigl(\mX, \bigoplus\limits_{i \in I} \mE_i[-n]\bigr) \cong
 \Hom_{\sA}\bigl(H^n(\mX), \bigoplus\limits_{i \in I} \mE_i\bigr)
 \cong
 \bigoplus\limits_{i \in I} \Hom_{\sA}\bigl(H^n(\mX), \bigoplus\limits_{i \in I} \mE_i\bigr).
 $$
 A result of Rentschler \cite{Rentschler} allows to conclude  that $H^n(\mX)$ is Noetherian.
\end{proof}

\noindent
The following theorem is the main result of our article.

\begin{theorem}\label{T:Main}
Let $\sX$  be a reduced rational projective curve with
nodal or cuspidal singularities, $\kA$ be its Auslander sheaf of orders and
$\kH^\bullet = \kS[-1] \oplus \kP(-1) \oplus \kP$ be the rigid complex from
Proposition \ref{T:crucial}. Then $\kH^\bullet$
is a tilting complex in the derived category $D^b\bigl(\Coh(\kA)\bigr)$.
\end{theorem}

\begin{proof}
By Proposition \ref{T:crucial}, the complex
$\kH^\bullet$ is rigid. In order to apply Theorem \ref{T:Keller}, we have
 to show that the right orthogonal
of $\kH^\bullet$ in the unbounded derived category $D\bigl(\Qcoh(\kA)\bigr)$ is zero.
Let $\sC = \cD(\kH^\bullet)$ be the smallest triangulated subcategory of  $D\bigl(\Qcoh(\kA)\bigr)$ containing  $\kH^\bullet$. Our goal is to show that
the right orthogonal of $\sC$ inside of $D\bigl(\Qcoh(\kA)\bigr)$ is zero.

\vspace{2mm}
\noindent
First observe that
 the Euler sequences
$$
0 \lar  \kO_{\mathbb{P}^1}(m-1) \lar  \kO_{\mathbb{P}^1}(m)^{\oplus 2}
\lar \kO_{\mathbb{P}^1}(m+1) \lar  0
$$
in the category of coherent sheaves on $\mathbb{P}^1$
induces the  short exact sequences
$$
0 \lar \kP_i(m-1) \lar  \kP_i(m)^{\oplus 2} \lar  \kP_i(m+1) \lar 0
$$
in the category $\Coh(\kA)$
for all $1 \le i \le n$ and $m \in \mathbb{Z}$. This implies
that all locally projective $\kA$-modules $\kP_i(m)$ belong to the
triangulated category $\sC$. In particular, the locally projective $\kA$--module $\kP(-2)$
belongs to the category $\sC$. The short exact sequence
$$
0 \lar \kP(-2) \lar \kF
\lar
\kS_1 \oplus \kS_2 \oplus \dots \oplus \kS_t \lar   0
$$
implies that the locally projective $\kA$--module $\kF$ belongs to $\sC$.

Consider the torsion $\kA$--module $\kT$, which is  the cokernel
of the canonical inclusion morphism $\kF \rightarrow \kP$. Since
$\kF$ and $\kP$ belong to $\sC$, it follows that  $\kT$ is an object of
 $\sC$, too.
Moreover, $\kT$ is  supported at the singular locus
of $\sX$, hence
$
\kT \cong \displaystyle \bigoplus\limits_{j=1}^t \kT_{p_j}.
 $
Let $\mT_j$ be the finite-length module over $\widehat{\kA}_{p_j}$ corresponding to
the torsion sheaf $\kT_{p_j}$.
We have:
\begin{itemize}
\item if $p_j$ is a node  then $\mT_j$ is given by
$$
\xymatrix
{
\kk  \ar@/^/[rr]^{1}  & &  \kk \ar@/^/[ll]^{0}
 \ar@/_/[rr]_{0}
 & &
\ar@/_/[ll]_{1} \kk}
$$
\item
if $p_j$ is a  cusp then $\mT_j$ is given by the quiver representation
$$
\xymatrix
{
\kk^2
\ar@(ul, dl)_{
\left(
\begin{smallmatrix}
0 & 1 \\
0 & 0
\end{smallmatrix}
\right)
} \ar@/^/[rr]^{\left(\begin{smallmatrix} 1 0  \end{smallmatrix}\right)}  & & \kk \ar@/^/[ll]^{
\left(
\begin{smallmatrix}
0 \\
0
\end{smallmatrix}
\right)
}
}
$$
\end{itemize}
If the point $p_j$ is nodal then for any $\lambda \in \kk^*$ we have a short exact sequence
$$
0 \lar  \mT_j \lar \mU_j(\lambda) \lar \mS_j \lar 0,
$$
where
$\mU_j(\lambda)$ is the module
$$
\xymatrix
{
\kk  \ar@/^/[rr]^{\left(\begin{smallmatrix} 0 \\ 1\end{smallmatrix}\right)}  & &
\kk^2 \ar@/^/[ll]^{\left(\begin{smallmatrix} 1 0 \end{smallmatrix}\right)}
 \ar@/_/[rr]_{\left(\begin{smallmatrix} 1 0 \end{smallmatrix}\right)}
 & &
\ar@/_/[ll]_{\left(\begin{smallmatrix} 0 \\ \lambda\end{smallmatrix}\right)} \kk}.
$$
Similarly, for a  cuspidal point $p_j$  the module $\mU_j(\lambda)$ is given by the representation
$$
\xymatrix
{
\kk^2
\ar@(ul, dl)_{
\left(
\begin{smallmatrix}
0 & 1 \\
0 & 0
\end{smallmatrix}
\right)
} \ar@/^/[rr]^{\left(\begin{smallmatrix} 0  &  0 \\ 1  & 0  \end{smallmatrix}\right)}  & & \kk^2 \ar@/^/[ll]^{
\left(
\begin{smallmatrix}
\lambda & 0  \\
1 & 0
\end{smallmatrix}
\right)
}
}
$$
for some $\lambda \in \kk$.
In particular, for any choice of $\lambda_1, \lambda_2, \dots, \lambda_t \in \kk^*$,
the torsion sheaf $\kU$ corresponding to the module $\bigoplus\limits_{j = 1}^t \mU_j(\lambda_j)$
belongs to the category $\sC$.

Let $\kL \in \Pic(\sX)$ be a line bundle of multidegree $(1, 1, \dots, 1)$. Then for any $m \ge 1$ there
exists an exact sequence
$$
0 \lar \kF \otimes_\kO \kL^{-m} \lar \kF \otimes_\kO \kL^{-m+1} \lar
\bigoplus\limits_{j = 1}^t \kU_j(\lambda_j) \lar 0.
$$
In particular, for all $m \ge 0$ the sheaf $\kF \otimes_\kO \kL^{-m}$ belongs to $\sC$.
Hence, for any $m \ge 0$ the locally free $\kA$--module
$\kF \otimes_\kO \kL^{-m}$ belongs to $\sC$, too.

Let $\kG$ be an arbitrary coherent $\kA$--module. Then it is also a coherent $\kO$-module.
Since $\kL$ is an ample line bundle on $\sX$, by a theorem of Serre
 there exists $m\ge 0$ such that the evaluation morphism
 $\Hom_\sX(\kL^{-m}, \kG) \stackrel{\kk}\otimes \kL^{-m}
 \stackrel{\mathsf{ev}}\lar \kG$ is surjective.
In particular, there exists $N \ge 0$  such that there exists an epimorphism of $\kO$--modules
$(\kL^{-m})^{N} \rightarrow \kG$. Since the functor $\kA \otimes_\kO \,-\,$ is right exact,
we have epimorphisms
$
\xymatrix{
\kA \otimes_\kO (\kL^{-m})^{N} \ar@{->>}[r]^-{1 \otimes \mathsf{ev}} &
\kA \otimes_\kO \kG  \ar@{->>}[r]^-{\mathsf{can}} &   \kG.
}
$

Let $\kG^\bullet = \bigl(\dots \rightarrow \kG^{n-1} \xrightarrow{\delta^{n-1}}
\kG^{n} \xrightarrow{\delta^{n}} \kG^{n+1} \rightarrow \dots \bigr)
$ be a complex from the category $D\bigl(\Qcoh(\kA)\bigr)$. If $\kG^\bullet \not\cong 0$ then
there exists $n \in \ZZ$ such that $\kH^n(\kG^\bullet) \ne 0$. Let $\kK_n = \ker(\delta^n)$.
Since any quasi-coherent $\kA$--module is a direct limit of its coherent submodules, there
exists a \emph{coherent} submodule $\kN_n$ of $\kK_n$ such that it is not a subobject
of the sheaf $\im(\delta^{n-1}) \subseteq \kK_n$.

Consider the morphism $\bigl(\kA \otimes \kL^{-m}\bigr)^N[-n]
\xrightarrow{f} \kG^\bullet$ in the derived category
$D\bigl(\Qcoh(\kA)\bigr)$ defined as the composition
$$
\bigl(\kA \otimes \kL^{-m}\bigr)^N[-n] \lar
\kN_n[-n] \lar \kK_n[-n] \lar \kG^\bullet.
$$
Since  $\kH^n(f) \ne 0$, the morphism $f$ is non-zero in $D^b\bigl(\Qcoh(\kA)\bigr)$, too. This shows that
$\sC^\perp = 0$, hence  $\kH^\bullet$ is a tilting complex.
\end{proof}

\begin{corollary}\label{C:embedofcategor}
Theorem \ref{T:Main}, Theorem \ref{T:Keller} and Proposition \ref{P:Murfet} imply that
there exists an equivalence of triangulated categories
$
\TT: D\bigl(\Qcoh(\sX)\bigr) \lar D\bigl(\mathsf{Mod}-\Gamma_\sX\bigr)
$
inducing equivalences of triangulated categories
$$
D^{-}\bigl(\Coh(\kA)\bigr) \lar
D^{-}\bigl(\mod-\Gamma_\sX\bigr) \,\, \mbox{and} \, \,
D^{b}\bigl(\Coh(\kA)\bigr) \lar
D^{b}\bigl(\mod-\Gamma_\sX\bigr).
$$
In particular, we have exact  fully faithful functors
$$
D^{-}\bigl(\Coh(\sX)\bigr) \lar
D^{-}\bigl(\mod-\Gamma_\sX\bigr) \,\, \mbox{and} \, \,
\Perf(\sX)  \lar
D^{b}\bigl(\mod-\Gamma_\sX\bigr).
$$
\end{corollary}

\subsection{Description of the tilted algebra}
Our next goal is to describe the tilted algebra
$\End_{D^b(\kA)}(\kH^\bullet)$ as the path algebra of some quiver with
relations.
Recall our notation.
Let $\sX$ be a rational projective curve with only nodes and cusps
as singularities, $\pi: \widetilde\sX \to \sX$
its normalization,
$\widetilde\sX =  \bigcup_{i=1}^n \widetilde{\sX_i}$, where
all $\widetilde\sX_i \cong {\mathbb P}^1$. Let
$$\mathsf{Sing}(\sX) = \bigl\{p_1, p_2, \dots,p_r, p_{r+1},
\dots,  p_{r+s}\bigr\}$$ be the singular locus of $\sX$,
 where $p_1, \dots, p_r$ are nodes and
$p_{r+1}, \dots, p_{r+s}$ are cusps.
Choose homogeneous coordinates $(u_i: v_i)$ on each irreducible component
$\widetilde{\sX_i}$ and  for any pair $(i,j)$ such that
$1 \le i \le n$ and $1 \le j \le r+s$ consider the set of points
$$
\pi^{-1}(p_j) \cap \widetilde{\sX_i} =
\bigl\{
q_{ij}^{(k)}  = (\alpha_{ij}^{(k)}:\beta_{ij}^{(k)}) \, \big|\,
1 \le k \le m_{ij}
\bigr\},
$$
where $m_{ij} = 0,1$ or $2$.
We additionally
assume the  coordinates are chosen in such a way that
$(\alpha_{ij}^{(k)}:\beta_{ij}^{(k)}) \ne (1: 0)$ for all
indices $i,j,k$ such that $p_j$ is a cusp.

\begin{definition}\label{D:tiltedalg}
The algebra $\Gamma_\sX$ attached to a rational projective curve $\sX$
with nodal or cuspidal singularities, is the path algebra of  the
following quiver with relations:
\begin{itemize}
\item It has   $2n + r +s$ vertices: for each index $1 \le i \le n$ we
have two points $a_i$ and $b_i$ and  for each index $1 \le j \le r+s$ we have
one point $c_j$.
\item The arrows of $\Gamma_\sX$ are as follows.
\begin{itemize}
\item For any index $1 \le i \le n$ we have two arrows $u_i, v_i: a_i \to b_i$;
\item For any index $1 \le j \le r$ (nodal point)
we have $m_{ij}$ arrows $w_{ij}^{(k)}: c_j \to a_i$;
\item For any index $r+1 \le j \le r+s$ (cuspidal point) and the
  unique index $i$ with $m_{ij}=1$ we have two arrows
$w_{ij}, w_{ij}' : c_j \to a_i$
\end{itemize}
\item The relations are as follows: for any $1 \le i \le n$,
$1 \le j \le r$ (nodal point) and $1 \le k \le m_{ij}$ we have:
$$
\bigl(\beta_{ij}^{(k)}u_i - \alpha_{ij}^{(k)}v_i\bigr)w_{ij}^{(k)} = 0
$$
and for any $1 \le i \le n$,
$r+1  \le j \le r +s$, $m_{ij}=1$ (cuspidal  point) we have:
$$
\bigl(\beta_{ij}^{(1)}u_i - \alpha_{ij}^{(1)}v_i\bigr)w_{ij}' = 0,
\quad \bigl(\beta_{ij}^{(1)}u_i - \alpha_{ij}^{(1)}v_i\bigr)w_{ij} = v_i w_{ij}'.
$$
\end{itemize}
\end{definition}

\begin{remark} The algebra $\Gamma_\sX$ is exactly the algebra defined
in \cite[Appendix A]{DGVB}.
Since  all paths in  $\Gamma_\sX$ have the length at most two, we have:
$\gldim(\Gamma_\sX) = 2$.
\end{remark}

\begin{proposition}\label{P:Main}
 The endomorphism
algebra of the tilting complex $\kH^\bullet$  from Theorem \ref{T:Main} is isomorphic to the algebra
$\Gamma_\sX$ introduced in Definition \ref{D:tiltedalg}.
\end{proposition}

\begin{proof}
A choice of  homogeneous coordinates
 $(u_i: v_i)$ on $\widetilde{\sX}_i$ yields a pair of distinguished sections
 $z_i^0, z_i^\infty \in H^0\bigl(\kO_{\widetilde{X}_i}(1)\bigr)$, where
 $z_i^0(0:1) = 0$ and $z_i^\infty(1:0) = 0$.
 They correspond to
a pair  of distinguished morphisms  $
u_i, v_i \in  \Hom_\kA\bigl(\kP_i(-1), \kP_i\bigr), \quad 1 \le i \le n
$
and  form  a basis of this morphism space.
In the course of the proof of Proposition  \ref{T:crucial} we have seen
that the only non-trivial contributions to $\End_{D^b(\Coh(\kA))}(\kH^\bullet)$ come
from:
$$
\Hom_\kA\bigl(\kP_i(-1), \kP_i\bigr) \cong
\Hom_{\mathbb{P}^1}\bigl(\kO_{\mathbb{P}^1}(-1), \kO_{\mathbb{P}^1}\bigr)
= \kk^2
$$
and
$$
\Ext^1_\kA\bigl(\kS_j, \kP_i(-1)\bigr) \cong
H^0\Bigl({\mathcal Ext}^1_\kA\bigl(\kS_j, \kP_i(-1)\bigr)\Bigr) \cong \kk^2 \cong
\Ext^1_\kA\bigl(\kS_j, \kP_i\bigr).
$$
Since the spaces $\Ext^1_\kA\bigl(\kS_j, \kP_i(-1)\bigr)$ can be computed locally,
we carry out calculations  over the complete ring $\widehat{\kA}_{p_j}$, following
 the notations of Remark \ref{R:Auslanderalg}.

\medskip
\noindent
\underline{1-st case}. Assume $p_j \in \mathsf{Sing}(\sX)$ is nodal and its both preimages
$$\pi^{-1}(p_j) = \bigl\{q_{ij}^{(1)}, q_{ij}^{(2)}\bigr\} =
\bigl\{(\alpha_{ij}^{(1)}: \beta_{ij}^{(1)}),
(\alpha_{ij}^{(2)}: \beta_{ij}^{(2)})\bigr\}$$
belong to the same irreducible component  $\widetilde{\sX_i}$. Recall that in this case,
the algebra $\widehat{\kA_{p_j}}$ is given by the completion of the following
quiver with relations:
$$
\xymatrix
{
\stackrel{1}\bullet \ar@/^/[rr]^{a^{-}_j}  & &  \stackrel{2}\bullet \ar@/^/[ll]^{a^{+}_j}
 \ar@/_/[rr]_{b^{+}_j}
 & &
\ar@/_/[ll]_{b^{-}_j} \stackrel{3}\bullet}  \qquad b^{+}_j a^{-}_j = 0,
\quad  a^{+}_j b^{-}_j = 0.
$$
Moreover, we have:  $\widehat{(\kP_i)}_{p_j} \cong \mP_j^{(1)} \oplus \mP_j^{(3)}$, $\widehat{\kF}_{p_j} \cong \mP_j^{(2)}$.
Recall that the simple module
$\mS_j = \mS_j^{(2)}$
has a projective resolution
$$
0 \lar  \mP_j^{(1)} \oplus \mP_j^{(3)} \xrightarrow{(a^{+}_j b^{+}_j)}
\mP_j^{(2)} \lar \mS_j \lar  0.
$$
Hence,
$
\Ext^1_{\widehat{\kA}_{p_j}}\bigl(\mS_j, \mP_j^{(1)} \oplus \mP_j^{(3)}\bigr) \cong
\kk^2 = \bigl\langle w_{ij}^{(1)}, w_{ij}^{(2)}\bigr\rangle,
$
where $w_{ij}^{(1)}$ is given by the following morphism in the homotopy category:
$$
\xymatrix
{ 0 \ar[r] & \mP_j^{(1)} \oplus \mP_j^{(3)} \ar[r] \ar[d]_{
\left(
\begin{smallmatrix}
1 & 0 \\
0 & 0
\end{smallmatrix}
\right)
} & \mP_j^{(2)} \\
  0 \ar[r] & \mP_j^{(1)} \oplus \mP_j^{(3)} \ar[r] & 0
}
$$
and the morphism $w_{ij}^{(2)}$ is defined in a similar way.
This implies that any morphism from $\kS_j[-1]$ to $\kP_i$
in the derived category $D^b\bigl(\Coh(\kA)\bigr)$
factors through $\kP_i(-1)$.
Next, note that the section $\beta_{ij}^{(k)} u_i - \alpha_{ij}^{(k)} v_i$ vanishes
only at the point $q_{ij}^{(k)} = (\alpha_{ij}^{(k)}: \beta_{ij}^{(k)})$, $k = 1,2$.  Hence,
we have the equalities:
\begin{equation}\label{E:node}
(\beta_{ij}^{(1)} u_i - \alpha_{ij}^{(1)} v_i) w_{ij}^{(1)} = 0, \quad
(\beta_{ij}^{(2)} u_i - \alpha_{ij}^{(2)} v_i) w_{ij}^{(2)} = 0
\end{equation}
in the morphism space $\Hom_{(\kA)}\bigl(\kS_j[-1], \kP_i\bigr)$.
Moreover, the morphisms
$(\beta_{ij}^{(1)} u_i - \alpha_{ij}^{(1)} v_i) w_{ij}^{(2)}$
and $(\beta_{ij}^{(2)} u_i - \alpha_{ij}^{(2)} v_i) w_{ij}^{(1)}$ are linearly independent.
Since $\Ext^1_{\widehat{\kA}_{p_j}}\bigl(\mS_j, \mP_j^{(1)} \oplus \mP_j^{(3)}\bigr) \cong
\kk^2$, there are no other relations between $w_{ij}^{(k)}$, $u_i$ and $v_i$ but those described
in  (\ref{E:node}).

\medskip
\noindent
\underline{2-nd case}. Assume the point $p_j$ is cuspidal and
$\pi^{-1}(p_j) = q_{ij} = (\alpha_{ij}: \beta_{ij}) \in \widetilde{\sX}_i$.
 In this case, $\widehat{\kA}_{p_j}$ is isomorphic
to  the completion of the path algebra
$$
\xymatrix
{
\stackrel{1}\bullet \ar@(ul, dl)_{a_j} \ar@/^/[rr]^{b^{+}_j}  & &
\stackrel{2}\bullet \ar@/^/[ll]^{b^{-}_j}
} \qquad a_j^2 = b^{-}_j b^{+}_j.
$$
We have:
$\bigl(\kP_i(-1)\bigr)_{p_j} \cong (\kP_i)_{p_j} \cong \mP_j^{(1)}$
and
$
0 \to \mP_j^{(1)}    \stackrel{b^{-}_j}\lar \mP_j^{(2)} \to \mS_j \to 0
$
is a projective resolution of the \emph{rigid} simple module $\mS_j$.
Hence,
$
\Ext^1_\kA(\kS_j, \kP_i) \cong  H^0\bigr({\mathcal Ext}^1_\kA(\kS_j, \kP_i)\bigr) = \kk^2.
$
Since the homogeneous coordinates $(u_i: v_i)$ are chosen in such a way
that $q_{ij} \ne (1: 0)$, the morphisms $v_i$ and
$\beta_{ij} u_i - \alpha_{ij} v_i: \kP_i(-1) \rightarrow   \kP_i$
are linearly independent. Let
$w_{ij}':= (\beta_{ij} u_i - \alpha_{ij} v_i)_{p_j}$ and
$w_{ij}:= (v_i)_{p_j}$ be the induced $\widehat{\kA}_{p_j}$--linear morphisms
of the projective module $\mP^{(1)}_j$. Denote by the same letters the induced morphisms
of complexes from $\bigl(\mP_j^{(1)}    \stackrel{b^{+}_j}\lar \mP_j^{(2)}\bigr)[-1]$ to
$\mP_j^{(1)}$
in the homotopy category
of projective $\widehat{\kA}_{p_j}$--modules. Then
$w_{ij}$ and $w_{ij}'$ form a basis of $\Ext_\kA^1\bigl(\kS_j, \kP_i(-1)\bigr) = \kk^2$ and
 we obtain the relations
$$
\left\{
\begin{array}{ccc}
(\beta_{ij} u_i - \alpha_{ij} v_i) w_{ij}' & = & 0 \\
(\beta_{ij} u_i - \alpha_{ij} v_i) w_{ij} & = & v_i w_{ij}'
\end{array}
\right.
$$
in the morphism space $\Ext^1_\kA(\kS_j, \kP_i)$. As in the Case 1, it follows that there
are no other relations between $u_i, v_i, w_{ij}$ and $w_{ij}'$.

\medskip
\noindent
\underline{3-rd case}.
The case when $p_j$ is nodal and its preimages belong to different components
of $\widetilde\sX$ is completely similar to the first case and is therefore left to the reader.
\end{proof}

\begin{example}
Let $\sX$ be an irreducible nodal rational projective curve of arithmetic genus two,
$p_1$ and $p_2$ its singular points, $\PP^1 \xrightarrow{\pi} \sX$ its normalization.
Assume that coordinates on $\PP^1$ are chosen in such a way that
$$\pi^{-1}(p_1) = \bigl\{0 = (0:1), \infty = (1:0)\bigr\} \quad  \mbox{and} \quad
\pi^{-1}(p_2) = \bigl\{(1:1), (\lambda: 1)\bigr\}$$ with $\lambda \in \kk\setminus \{0, 1\}$.
Then the algebra $\Gamma_\sX$ is the path algebra of the following quiver
$$
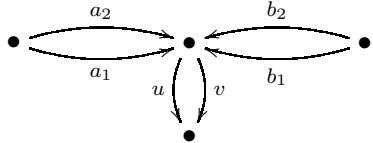
\xymatrix
{
\bullet \ar@/_/[rr]_{a_1} \ar@/^/[rr]^{a_2} &  & \bullet \ar@/^/[d]^v  \ar@/_/[d]_u & & \bullet \ar@/^/[ll]^{b_1} \ar@/_/[ll]_{b_2}  \\
                                 & & \bullet  & &
}
$$
subject to the relations $u a_1 = 0, \, v a_2 = 0, \, (u-v) b_1 = 0$ and $(u- \lambda v) b_2 = 0$.
It seems to be an interesting problem to study  compactified moduli spaces of vector bundles
on $\sX$ in terms of representations of the algebra $\Gamma_\sX$.
\end{example}

\subsection{Dimension of the derived category of a rational
 projective curve}
As a consequence of our approach, we  obtain an upper bound of the dimension
of the derived category of coherent sheaves of reduced rational projective curve
with only nodal or cuspidal singularities.

Let $\sC$ be an idempotent complete triangulated subcategory and $\sA$, $\sB$
 be its two idempotent complete full subcategories closed under shifts.
 Following Rouquier \cite{Rouquier}, we denote by $\sA * \sB$ the full subcategory
 of $\sC$, whose  objects are those objects  $\mX$  of  $\sC$ for which there exists
a distinguished triangle
$$
\mA \lar \mX \lar \mB\lar \mA[1]
$$
with $\mA \in \Ob(\sA)$ and $\mB \in \Ob(\sB)$. For an object $\mX \in \Ob(\sC)$
we denote by $\langle \mX \rangle$ the smallest full subcategory of $\sC$, closed under
taking shifts, direct sums and direct summands. Next, for any positive integer $n$
we define subcategories $\langle \mX\rangle_n$ by the following rule:
$$
\langle \mX\rangle_1 = \langle \mX\rangle, \quad \langle \mX\rangle_{n+1} =
\bigl\langle\langle \mX\rangle_1 * \langle \mX\rangle_n\bigr\rangle.
$$
An object $\mX \in \Ob(\sC)$ is a   \emph{strong generator}
if  $\langle \mX\rangle_n = \sC$ for some positive integer $n$.
Rouquier suggested the following definition of the dimension of a triangulated category
$\sC$:
$$
\dim(\sC) = \inf \left\{n \in \ZZ_+  \, \Big| \, \exists \,  \mX \, \in \Ob(\sC): \langle \mX\rangle_{n+1} = \sC \right\}.
$$
He has also proven  that the dimension of the derived category of coherent sheaves of a separated
scheme    $\sX$  of finite type over a perfect field $\kk$ is always finite, see \cite[Theorem 7.38]{Rouquier}.
Moreover, if $\sX$ is smooth of dimension $n$ then
$
n \le \dim\bigl(D^b\bigl(\Coh(\sX)\bigr)\bigr) \le 2n,
$
see \cite[Proposition 7.9 and Proposition 7.16]{Rouquier}.
By a recent result of Orlov \cite{Orlov}, for a \emph{smooth} projective curve $\sX$ over a field $\kk$ we have:
$$
\dim\bigl(D^b\bigl(\Coh(\sX)\bigr)\bigr) = 1.
$$
The case of the singular projective curves still remains open. However, our technique allows to deduce
the following result.

\begin{theorem}
Let $\sX$ be a reduced rational projective curve with only nodal or cuspidal singularities.
Let $S = \{p_1, p_2,\dots, p_t\}$ be the singular locus of $\sX$ and $\sX_1, \sX_2, \dots,
\sX_n$ be the irreducible components of $\sX$. Let $\kO_i = \kO_{\sX_i}$ be the structure sheaf
of $\sX_i$, $1 \le i \le n$ and $\kO_i(-1) = \kO_i \otimes \kO_\sX(-q_i)$, where
$q_i \in \sX_i$ is a smooth point. Consider the coherent sheaf
$$
\kG = \bigoplus\limits_{i=1}^n \Bigl(\kO_i(-1) \oplus \kO_i\Bigr) \oplus \bigoplus\limits_{j=1}^t \kk_{p_j}.
$$
Then we have: $\langle \kG\rangle_3 = D^b\bigl(\Coh(\sX)\bigr)$. In particular,
$
\dim\bigl(D^b\bigl(\Coh(\sX)\bigr)\bigr) \le 2.
$
\end{theorem}

\begin{proof}
Let $\kA$ be the Auslander sheaf of  $\sX$. By Theorem \ref{T:Main},  the derived category
$
D^b\bigl(\Coh(\kA)\bigr)
$
is equivalent to $D^b(\mod-\Gamma_\sX)$. Moreover, $\gldim(\Gamma_\sX) = 2$ and
the equivalence $\TT$ maps the tilting complex $\kH^\bullet$ to the regular module
$\Gamma_\sX$. By \cite[Lemma 7.1]{Rouquier} it is known that
$\langle \Gamma_\sX\rangle_3 = D^b(\mod-\Gamma_\sX)$. This implies that
$
\langle \kH^\bullet \rangle_3 = D^b\bigl(\Coh(\kA)\bigr).
$

Consider now the exact functor $\GG: \Coh(\kA) \rightarrow  \Coh(\sX)$.  By Theorem \ref{T:main},
the derived functor $\GG: D^b\bigl(\Coh(\kA)\bigr) \rightarrow  D^b\bigl(\Coh(\sX)\bigr)$
is essentially surjective, hence
$
\bigl\langle \GG(\kH^\bullet) \bigr\rangle_3 = D^b\bigl(\Coh(\sX)\bigr).
$
To conclude the proof, it remains to note that
$$
\GG(\kH^\bullet) \cong \bigoplus\limits_{i=1}^n \Bigl(\kO_i(-1) \oplus \kO_i\Bigr) \oplus \bigoplus\limits_{j=1}^t \kk_{p_j}[-1].
$$
\end{proof}

\section{Coherent sheaves on Kodaira cycles and gentle algebras}

In this section we discuss some corollaries from the results obtained
in the previous section, in the case of  Kodaira cycles of projective
lines. To deal with \emph{left modules}, we prefer to replace the tilted
algebra $\Gamma_X$ by its opposite $\Lambda_X=\Gamma_X^{\mathrm{op}}$.

\begin{proposition}
Let $\sE = \sE_n$ be a Kodaira cycle of $n$ projective lines (in the
case $n = 1$ it is an irreducible plane nodal cubic curve),  $\kA$ be the
Auslander sheaf  and
$\Lambda = \Lambda_\sE$ be the opposite algebra of
the corresponding tilted algebra.
 Then we have:
\begin{enumerate}
\item The algebra $\Lambda_\sE$ is gentle, see
\cite{AssemSkowr}.
\item The categories  $\Perf(\sE)$  and $\Coh(\sE)$ are   tame in the
  ``pragmatic sense''
\footnote{A more precise result implying the tameness in a ``strict
  sense'' was obtained in \cite{Duke}.}.
\end{enumerate}
\end{proposition}

\begin{proof} The fact that the algebra $\Lambda$ is gentle, follows
from Proposition \ref{P:Main} and the definition of the gentle algebras.
Moreover, the gentle algebras are derived-tame, see  \cite{PogSkowr} and \cite{Ringelrep}.
Since we have  fully faithful functors
$\Perf(\sE)
\rightarrow  D^b\bigl(\Coh(\kA)\bigr) \xrightarrow{\sim}  D^b\bigl(\Lambda-\mod\bigr),
$
and
$
\Coh(\sE) \rightarrow \Coh(\kA) \rightarrow  D^b\bigl(\Coh(\kA)\bigr) \xrightarrow{\sim}  D^b\bigl(\Lambda-\mod\bigr),
$
the categories $\Perf(\sE)$ and $\Coh(\sE)$ are equivalent to full subcategories of
a representation-tame category $D^b\bigl(\Lambda-\mod\bigr)$. This precisely means they
are pragmatic-tame.
\end{proof}

\begin{example}
Let $\sE = \sE_2$ be a Kodaira cycle of two projective lines. Then the algebra $\Lambda_\sE$
is the  path algebra of  the following quiver
$$
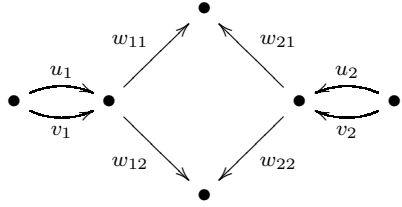
\xymatrix
{
 & & \bullet & & \\
\bullet \ar@/^/[r]^{u_1} \ar@/_/[r]_{v_1}
 & \bullet \ar[ru]^{w_{11}} \ar[rd]_{w_{12}} & &
 \bullet \ar[lu]_{w_{21}} \ar[ld]^{w_{22}}
& \bullet \ar@/_/[l]_{u_2} \ar@/^/[l]^{v_2} \\
& & \bullet & &
}
$$
subject to the relations:
$
w_{11} v_1 = 0, \, w_{22} v_2 = 0, \, w_{12} u_1 = 0$ and $w_{21} u_2 = 0.
$
In particular, there exists  a fully faithful functor
$
\Perf(\sE) \rightarrow   D^b\bigl(\Lambda_\sE-\mod\bigr).
$
\end{example}

\noindent
Let $\Lambda$ be a finite-dimensional algebra over a field $\kk$. Then the Nakayama functor
$$\nu:= \DD \Hom_\Lambda(\,-\,,\Lambda): \Lambda-\mod \lar  \Lambda-\mod$$ is right exact.
Moreover, if $\Lambda = \kk\vec{Q}/\rho$ is the path algebra of a finite quiver with relations  then
$\nu(\mP_i) = \mI_i$, where $\mP_i$ and $\mI_i$ are the indecomposable projective and
injective modules corresponding to the vertex $i \in Q_0$.
If  $\gldim(\Lambda) < \infty$ then
by a result of Happel \cite{Happel}, the derived functor
$$
\Ss:= \mathbb{L}\nu: D^b(\Lambda-\mod) \lar D^b(\Lambda-\mod)
$$
is the Serre functor of the category $D^b(\Lambda-\mod)$ and the
functor $\tau_\Lambda = \Ss[-1]$ is the Auslander-Reiten translate
in $D^b(\Lambda-\mod)$ (see also \cite{ReitenvandenBergh}).

\begin{corollary}\label{C:appltogentalg}
Let $\sE$ be a Kodaira cycle of projective lines, $\kA$ be its  Auslander sheaf
and $\Lambda$ be the opposite algebra of the corresponding tilted algebra.
Consider the category $\Band(\Lambda)$,
which is  the full subcategory of $D^b(\Lambda-\mod)$ whose objects are  the complexes
$\mP^\bullet$ such that $\tau_\Lambda(\mP^\bullet) \cong \mP^\bullet$. Then $\Band(\Lambda)$ is triangulated and
idempotent complete. Moreover, it is triangle equivalent to the category of perfect complexes
$\Perf(\sE)$.
\end{corollary}

\begin{proof}
Since $\TT: D^b\bigl(\Coh(\kA)\bigr) \rightarrow  D^b(\Lambda-\mod)$ is an equivalence of categories,
we have an isomorphism of functors $\TT \circ \tau_\kA \cong \tau_{\Lambda} \circ \TT$. This implies
that the category $\Band(\Lambda)$ is equivalent to the category $\cD_\theta$ introduced
in Corollary \ref{C:transformtheta}. By Theorem \ref{T:tauperiod} and Corollary \ref{C:tauperiod}
we get that $\Band(\Lambda)$ is equivalent to $\Perf(\sE)$. In particular, it is
idempotent-complete.
\end{proof}

\begin{remark}
Corollary \ref{C:appltogentalg}  implies the following result on the shape of the Auslander--Reiten quiver
of the gentle algebra $\Lambda$, attached to a cycle of projective lines
 $\sE$. Let $\mP^\bullet$ be an indecomposable object of
$D^b(\Lambda-\mod)$ such that $\tau_\Lambda^n(\mP^\bullet) \cong \mP^\bullet$. Then
$n = 1$ and $\mP^\bullet$ is an object of $\Band(\Lambda)$. In other words, the Auslander-Reiten
quiver of $D^b(\Lambda-\mod)$ does not contain tubes of length bigger than one.

Note that all indecomposable  \emph{band complexes} of
$D^b(\Lambda-\mod)$ in the sense of the work of Butler and Ringel
\cite{ButlerRingel} (see also \cite{Ringelrep}) are
$\tau$--periodic. However, their
definition of bands and strings of the derived category of a gentle algebra  is purely combinatorial.
In particular, in certain cases algebra automorphisms map
band complexes to string complexes. Moreover, it is not completely clear that
their notion of bands and strings coincide with the corresponding notion of bands and strings
 used in our previous paper \cite{Toronto}.

 In the particular situation of a gentle algebra $\Lambda$ which is the tilted algebra attached
 to a Kodaira cycle of projective lines, we say that an indecomposable
 object of $D^b(\Lambda-\mod)$ is a \emph{band}
 if and only if it belongs to $\Band(\Lambda)$.
 The indecomposable objects  of $D^b(\Lambda-\mod)$ which are not bands will
 be called \emph{strings}.
 The description of indecomposable objects in $D^b(\Lambda-\mod)$ obtained in
\cite{PogSkowr}, \cite{Ringelrep} and  \cite{Toronto}
implies  that strings do not have continuous moduli and are
classified by discrete parameters. The same concerns the indecomposable objects
of the derived category $D^b\bigl(\Coh(\sE)\bigr)$ which do not belong to $\Perf(\sE)$, see
\cite{Duke}.

The interplay between various categories occurring in our construction
 can be explained by the following diagram:
$$
\xymatrix{
D^-\bigl(\Coh(\sE)\bigr) \ar[r]^\FF & D^-\bigl(\Coh(\kA_\sE)\bigr) \ar[r]^\TT & D^-\bigl(\Lambda_\sE-\mod\bigr) \\
D^b\bigl(\Coh(\sE)\bigr)  \ar@{^{(}->}[u] \ar@{^{(}->}[ru] & D^b\bigl(\Coh(\kA_\sE)\bigr) \ar[r]^\TT
\ar@{_{(}->}[u] & D^b\bigl(\Lambda_\sE-\mod\bigr) \ar@{_{(}->}[u] \\
\Perf(\sE) \ar[r]^\sim \ar@{^{(}->}[u] & \cD_\theta \ar[r]^\sim \ar@{_{(}->}[u] & \Band(\Gamma_\sE).  \ar@{_{(}->}[u]
}
$$
Summing up, we get  that the
triangulated category $\Perf(\sE)$ is a full subcategory of two different
derived categories.
From one side, it is a subcategory of the derived category of coherent sheaves
$\Coh(\sE)$, whose global dimension is \emph{infinity}. From another side, it
is a subcategory of the derived category of representations of the algebra
$\Gamma_\sE$, whose global dimension is \emph{two}. The ``complement'' of $\Perf(\sE)$ in both categories
is ``small'': it consists of direct sums of complexes described by discrete parameters.
\end{remark}

\section{Tilting exercises on a Weierstra\ss{} nodal cubic curve}

Let $\sE \subseteq \mathbb{P}^2$
be a singular Weierstra\ss{} cubic curve over an algebraically closed field
$\kk$ given by the equation
$zy^2 = x^3 + x^2z$.   By Corollary \ref{C:embedofcategor},  there exists  a fully faithful functor
$$
\Perf(\sE)  \stackrel{\mathbb{F}}\lar  D^b\bigl(\Coh(\kA)\bigr)
\stackrel{\TT}\lar D^b\bigl(\Lambda-\mod\bigr),
$$
where $\kA = \kA_\sE$ is  the  Auslander   sheaf of orders and $\Lambda = \Lambda_\sE$ is the path algebra
of the following quiver with relations:
$$
\xymatrix
{
1  \ar@/^/[rr]^{a} \ar@/_/[rr]_{c} & & 2 \ar@/^/[rr]^{b} \ar@/_/[rr]_{d} & & 3
}
\qquad ba = 0 = dc.
$$
The algebra $\Lambda$ is interesting from various perspectives. First of all, it is gentle, hence
derived-tame. Next, by  a work of Seidel \cite[Section 3]{Seidel}, it is related with the directed
Fukaya category of a certain Lefschetz pencil.

Our next goal is to compute the complexes in $D^b\bigl(\Lambda-\mod\bigr)$ corresponding to the
images of
certain  perfect coherent sheaves on $\sE$ under the functor $\TT \circ \FF$.

Let $\pi: \PP^1 \rightarrow \sE$ be the normalization of $\sE$ and $s = (0:0:1)  \in \sE$ be the singular point.
 Choose coordinates
on $\PP^1$ in such a way that $\pi^{-1}(s) = \{0, \, \infty\}$, where $0 = (0:1)$ and $\infty = (1:0)$.
This choice yields two distinguished sections $z_0, z_\infty \in \Hom_{\PP^1}\bigl(\kO_{\PP^1}(-1), \kO_{\PP^1}\bigr)
= H^0\bigl(\kO_{\PP^1}(1)\bigr)$ such that $z_0(0) = 0$ and $z_\infty(\infty) = 0$.  Recall that
$$
\rA := \widehat{\kA}_s \cong
\left(
\begin{array}{cc}
\kk\llbracket u \rrbracket \times \kk\llbracket v \rrbracket & (u, v) \\
\kk\llbracket u \rrbracket \times \kk\llbracket v \rrbracket & \kk\llbracket u,v \rrbracket/(uv)
\end{array}
\right)
$$
is isomorphic to the  completion of the path algebra of the  following quiver with relations:
$$
\xymatrix
{
\stackrel{\alpha}\bullet \ar@/^/[rr]^{a_{-}}  & &  \stackrel{\beta}\circ \ar@/^/[ll]^{a_{+}}
 \ar@/_/[rr]_{b_{+}}
 & &
\ar@/_/[ll]_{b_{-}} \stackrel{\gamma}\bullet}  \qquad b_{+} a_{-} = 0, \quad  a_{+} b_{-} = 0.
$$
Let
$$
\mP_\alpha =
\left(
\begin{array}{c}
\kk\llbracket u \rrbracket \\
\kk\llbracket u \rrbracket
\end{array}
\right), \,
\mP_\gamma =
\left(
\begin{array}{c}
\kk\llbracket v \rrbracket \\
\kk\llbracket v \rrbracket
\end{array}
\right) \, \, \mbox{and} \, \,
\mP_\beta = \mF  =
\left(
\begin{array}{c}
(u, v) \\
\kk\llbracket u,v \rrbracket/(uv)
\end{array}
\right)
$$
be the indecomposable projective $\rA$--modules. We distinguish  two locally projective
$\kA$--modules
$$
\kP =
\left(
\begin{array}{c}
\widetilde\kO \\
\widetilde\kO
\end{array}
\right)  \, \, \mbox{and} \, \,
\kF =
\left(
\begin{array}{c}
\kI \\
\kO
\end{array}
\right)
$$
and for any $n \in \ZZ$ we have: $\widehat{\kP(n)_s} \cong  \mP_\alpha \oplus \mP_\gamma$, whereas
$\widehat{\kF_s} \cong \mF$. By Lemma \ref{L:keyforcomp},  there are   the following canonical isomorphisms:
$$
H^0\bigl(\kO_{\PP^1}(1)\bigr) \cong H^0\Bigl(\kS\bigl(\kP(-1), \kP\bigr)\Bigr) \cong
\Hom_\kA\bigl(\kP(-1), \kP\bigr).
$$
Hence, we get two distinguished elements in $\Hom_\kA\bigl(\kP(-1), \kP\bigr)$, which will
be denoted by $z_0$ and $z_\infty$. They are  characterized by the property that there exists
isomorphisms $t_1: \widehat{\kP(-1)_s} \rightarrow  \mP_\alpha \oplus \mP_\gamma$ and
$t_2: \widehat{\kP_s} \rightarrow   \mP_\alpha \oplus \mP_\gamma$ making the following diagrams
commutative:
$$
\xymatrix
{
\widehat{\kP(-1)_s} \ar[r]^{z_0} \ar[d]_{t_1} & \widehat{\kP_s}  \ar[d]^{t_2} \\
\mP_\alpha \oplus \mP_\gamma
\ar[r]^{\left(\begin{smallmatrix} u & 0 \\ 0 & \mathsf{id} \end{smallmatrix}\right)}& \mP_\alpha \oplus \mP_\gamma
}
\quad
\xymatrix
{
\widehat{\kP(-1)_s} \ar[r]^{z_\infty} \ar[d]_{t_1} & \widehat{\kP_s}  \ar[d]^{t_2} \\
\mP_\alpha \oplus \mP_\gamma
\ar[r]^{\left(\begin{smallmatrix} \mathsf{id} & 0 \\ 0 & v \end{smallmatrix}\right)}& \mP_\alpha \oplus \mP_\gamma.
}
$$
Let $\kS_\beta$ be the torsion $\kA$--module supported at the singular point $s \in \sE$ and corresponding
to the simple $\rA$--module $\mS_\beta$. Then for any $n \in \ZZ$  the canonical map
$$
H^0\Bigl({\mathcal Ext}^1_\kA\bigl(\kS_\beta, \kP(n)\bigr)\Bigr) \lar
\Ext^1_\kA\bigl(\kS_\beta, \kP(n)\bigr)
$$
is an isomorphism. By Remark \ref{R:Auslanderalg}, the simple module
 $\mS_\beta$ has  the following projective resolution:
$$
0 \lar \mP_\alpha \oplus \mP_\gamma
\xrightarrow{\left(\begin{smallmatrix} u,\, v\end{smallmatrix}\right)} \mP_\beta \lar \mS_\beta \lar 0.
$$
Hence,  $\Ext^1_\rA(\mS_\beta, \mP_\alpha \oplus \mP_\gamma) = \kk^2 = \langle \xi, \eta\rangle$, where
$\xi$ and $\eta$ are induced by the  $A$--linear morphisms given by the matrices
$$
\xi = \left(
\begin{array}{cc}
1 & 0 \\
0 & 0
\end{array}
\right), \, \,
\eta = \left(
\begin{array}{cc}
0 & 0 \\
0 & 1
\end{array}
\right): \, \, \mP_\alpha \oplus \mP_\gamma \lar \mP_\alpha \oplus \mP_\gamma.
$$
In particular, $z_0 \xi = 0$ and $z_\infty \eta = 0$.
By the definition of the tilting equivalence
$\TT: D^b\bigl(\Coh(\kA)\bigr) \rightarrow D^b\bigl(\Lambda-\mod\bigr)$ given by the
 tilting complex $\kH^\bullet = \kS_\beta[-1] \oplus \kP(-1) \oplus \kP$,
 we have the following result.

\begin{lemma}
We have: $\TT(\kP) \cong \mP_1, \TT\bigl(\kP(-1)\bigr) = \mP_2$ and
$\TT(\kS_\beta) \cong \mP_3[1]$, where $\mP_i$ is the indecomposable projective $\Lambda$--module
corresponding to the vertex $i$, $1 \le i \le 3$.
\end{lemma}

\noindent
Our next goal is to compute the images of certain  finite length objects in
$\Coh(\kA)$.  Let $x = (\lambda: \mu) \in \PP^1 \setminus \{0, \infty\}$
 be an arbitrary point. Consider
the torsion $\kA$--module $\kT_x$ given by its  locally projective resolution
$$
0 \lar \kP(-1)
\xrightarrow{\mu z_0 - \lambda z_\infty} \kP \lar \kT_x \lar 0.
$$
Since  $x \not\in \{0, \infty\}$,   the sheaf
$\kT_x$ is supported at a smooth point of $\sE$ and $\Ext_\kA^i(\kS_\beta, \kT_x) = 0$ for
all $i \ge 0$.
Next, note that $\Hom_\kA\bigl(\kP(n), \kT_x\bigr) = \kk$ and $\Ext^i_\kA\bigl(\kP(n), \kT_x\bigr) = 0$
for all $i > 0$ and $n \in \ZZ$. This implies that $H^i\bigl(\TT(\kT_x)\bigr) = 0$ for
$i \ne 0$. In particular, the complex $\TT(\kT_x)$ is isomorphic in $D^b\bigl(\Lambda-\mod\bigr)$ to the stalk
complex
$$ \dots \lar 0 \lar   \mN_x \lar  0 \lar \dots
$$
where $\mN_x = H^0\bigl(\TT(\kT_x)\bigr)$.

Recall that for
an arbitrary representation $\mM$ of the quiver $\Lambda$ we have:
$\Hom_\Lambda(\mP_i, \mM) = \mM(i)$, where $\mM(i)$ is the dimension of $\mM$ at the vertex $i$,
$1 \le i \le 3$. This allows to compute the multi-dimension of the zero  cohomology
of $\TT(\kT_x)$:
 $\underline{\dim}(\mN_x) = \bigl(1, \, 1, \,0\bigr).$
Moreover, the right $\Gamma$--module $\mN_x$
endowed with the natural structure of an $\End_{D^b(\kA)}\bigl(\kH^\bullet)$--module
has the following projective resolution:
$$
0 \lar
\Hom_{D^b(\kA)}\bigl(\kH^\bullet, \kP(-1)\bigr)
\xrightarrow{(\mu z_0 - \lambda z_\infty)_*} \Hom_{D^b(\kA)}\bigl(\kH^\bullet, \kP\bigr) \lar \mN_x \lar 0.
$$
Interpreting it in the terms of quiver representations, we obtain the following result.

\begin{proposition}
For any $x = (\lambda: \mu) \in \PP^1 \setminus \{0, \infty\}$ we have: $\TT(\kT_x) = \mN_x[0]$, where
$\mN_x$ is the following representation of the algebra $\Lambda$:
$$
\mN_x := \mathsf{coker}\bigl(\mP_2 \xrightarrow{\mu a - \lambda c} \mP_1\bigr) \cong
\xymatrix
{
\kk  \ar@/^/[rr]^{\lambda} \ar@/_/[rr]_{\mu} & & \kk \ar@/^/[rr] \ar@/_/[rr] & & 0
}.
$$
\end{proposition}

\noindent
As the next step, we compute the images under $\TT$ of two other exceptional simple
$\kA$--modules.

\begin{proposition}
Let $\mS_\alpha$ and $\mS_\gamma$ be the simple $\rA$--modules corresponding  to
the vertices $\alpha$ and $\gamma$, and  $\kS_\alpha$ and $\kS_\gamma$ be
the corresponding torsion $\kA$--modules. Then we have:
$$\TT(\kS_\alpha) =
\xymatrix
{
\kk  \ar@/^/[rr]^{0} \ar@/_/[rr]_{1} & & \kk \ar@/^/[rr]^{1} \ar@/_/[rr]_{0} & & \kk
}
\,\,
\mbox{\rm and} \,\,
\TT(\kS_\gamma) =
\xymatrix
{
\kk  \ar@/^/[rr]^{1} \ar@/_/[rr]_{0} & & \kk \ar@/^/[rr]^{0} \ar@/_/[rr]_{1} & & \kk.
}
$$
\end{proposition}

\begin{proof}
First note that
$$\Ext^i_\kA(\kS_\beta, \kS_\alpha) = H^0\bigl({\mathcal Ext}^i_\kA(\kS_\beta, \kS_\alpha)\bigr) =
\left\{
\begin{array}{cl}
\kk & \mbox{if} \quad  i = 1, \\
0   & \mbox{otherwise}.
\end{array}
\right.
$$
In a similar way,  for any $n \in \ZZ$  we have:
$$
\Ext_\kA^i\bigl(\kP(n), \kS_\alpha\bigr) \cong H^0\Bigl({\mathcal Ext}_\kA^i\bigl(\kP(n), \kS_\alpha\bigr)\Bigr) =
\left\{
\begin{array}{cl}
\kk & \mbox{if} \quad  i = 0, \\
0   & \mbox{otherwise}.
\end{array}
\right.
$$
This implies that $\TT(\kS_\alpha)$ and $\TT(\kS_\gamma)$ are indecomposable stalk complexes.
Their zero cohomology are  $\Lambda$--modules $\mM_\alpha$ and $\mM_\gamma$,
whose multi-dimension is the vector  $(1, 1, 1).$ However, there are
precisely two \emph{indecomposable} $\Lambda$--modules with this  multi-dimension.
Since $\Hom_\kA(\kT_\infty, \kS_\alpha) = 0 = \Hom_\kA(\kT_0, \kS_\gamma)$ and
$\TT$ is an equivalence of categories, we have:
$\Hom_\Lambda(\mN_\infty, \mM_\alpha) = 0 = \Hom_\Lambda(\mN_0, \mM_\gamma).
$
This implies that $\mM_\alpha$ and $\mM_\gamma$ are given by the quiver representations as
stated above.
\end{proof}

\noindent
Finally, we shall compute the image of the Jacobian $\Pic^0(\sE)$ in the derived category
$D^b(\Lambda-\mod)$.

\begin{proposition}
The functor $\TT \circ \FF$ identifies the Jacobian  $\Pic^0(\sE) \cong \kk^*$ with  the following
family of complexes in the derived category $D^b(\Lambda-\mod)$
$$
\Bigl\{\mU^\bullet_\lambda\Bigr\}_{\lambda \in \kk^*} =
\left\{
\mP_3 \stackrel{\left(\begin{smallmatrix} \lambda b\\  d \end{smallmatrix}\right)}\lar
\underline{\mP_2 \oplus \mP_2} \stackrel{\begin{smallmatrix}(a, c) \end{smallmatrix}}\lar \mP_1
\right\}_{\lambda \in \kk^*},
$$
 where the underlined term of $\mU^\bullet_\lambda$
has degree zero.  Moreover,  for any
$\lambda \in \kk^*$ the complex $\mU^\bullet_\lambda$ is spherical in the sense
of Seidel and Thomas \cite{SeidelThomas} and we have: $H^0(\mU^\bullet_\lambda) = \mS_3$, $H^1(\mU^\bullet_\lambda) = \mS_1$
\end{proposition}

\begin{proof}
Since $\TT: D^b\bigl(\Coh(\kA)\bigr) \rightarrow D^b(\Lambda-\mod)$ is an equivalence
of categories, we have an isomorphism of functors $\TT \circ \tau_\kA \cong \tau_\Lambda \circ
\TT$. Next, we know that $\TT(\kP) = \mP_1$. This implies that
$$
\kR := \tau_\kA(\kP) =
\left(
\begin{array}{c}
\kI \\
\widetilde\kO
\end{array}
\right) \,\,
\mbox{and} \,\,
\TT(\kR) = \tau_\Lambda(\mP_1) = \mS_1[-1],
$$
where $\mS_1$ is the simple $\Lambda$--module corresponding to the vertex $1$
(note that $\mS_1$ is injective). Consider the torsion $\kA$--module $\kT$
supported at the singular point $s \in \sE$ and defined as
$$
\kT = \left(
\begin{array}{c}
0 \\
\widetilde\kO/\kI
\end{array}
\right)
=
\coker\left[
\left(
\begin{array}{c}
\kI \\
\kI
\end{array}
\right) \lar
\left(
\begin{array}{c}
\kI \\
\widetilde\kO
\end{array}
\right)
\right].
$$
Since $\Hom_\kA\bigl(\kP(n), \kS_\beta\bigr) = 0$ for all $n \in \ZZ$, the canonical
morphism
$$
\Hom_\kA(\kT, \kS_\beta) \lar \Hom_\kA(\kR, \kS_\beta)
$$
is an isomorphism. Another canonical morphism
$
\Hom_\kA(\kT, \kS_\beta) \lar \Hom_\sE(\widetilde\kO/\kI, \kO/\kI)
$
is an isomorphism as well, hence $\Hom_\kA(\kR, \kS_\beta) = \kk^2$. Since $(\widetilde\kO/\kI)_s \cong
\kk_s \times \kk_s$ and
$(\kO/\kI)_s = \kk_s$ are rings, the vector space $\Hom_\sE(\widetilde\kO/\kI, \kk_s)$ has two
\emph{distinguished} basis elements $\bar{w}_0$ and $\bar{w}_\infty$, which correspond to the non-trivial
idempotents of the ring $(\widetilde\kO/\kI)_s$. Let $w_0$ and $w_\infty$ be the corresponding elements
of $\Hom_\kA(\kR, \kS_\beta)$. For any $(\lambda,  \, \mu) \in \kk^2 \setminus \{0, 0\}$
consider the short exact sequence
$$
0 \lar \kX  \lar
\kR \stackrel{w}\lar  \kS_\beta \lar 0,
$$
where $w = \mu w_0 - \lambda w_\infty$. The $\kA$--module $\kX$ only depends
on the ratio $x= (\lambda: \mu) \in \PP^1$.

We claim that $\kX$ is a locally projective
$\kA$--module precisely when $x \not\in \{0, \infty \}$. Moreover, in this case
we have: $\kX_s \cong \mF$. Indeed, $\kX$ is locally projective at all smooth
points of $\sE$. Let $\bar{w} = \mu \bar{w}_0 - \lambda \bar{w}_\infty: \kT \rightarrow \kS_\beta$. Then
$\kX_s$ is isomorphic to the middle term of the following short exact sequence:
$$
0 \lar
\left(
\begin{array}{c}
\mI \\
\mI
\end{array}
\right) \lar
\left(
\begin{array}{c}
\mI \\
\mI + (\lambda, \mu) \rO
\end{array}
\right) \lar \ker(\bar{w}) \lar 0,
$$
where we view $(\lambda, \mu) \in \kk \times \kk$ as an element of
the normalization $\kk\llbracket u\rrbracket \times \kk\llbracket v\rrbracket$. Other way around, it is
not difficult to show that for any line bundle  $\kL_x \in \Pic^0(\sE) \cong \kk^*$ the locally
projective $\kA$--module $\kX := \kF \otimes_\sE \kL_x$ fits into  a short exact sequence
$$
0 \lar \kX \lar \kR \lar \kS_\beta \lar 0.
$$
Summing up, for all non-zero morphisms $w = \mu w_0 - \lambda w_\infty \in \Hom_\kA(\kR, \kS_\beta)$, where
$x = (\lambda: \mu) \notin \{0, \infty\}$, the mapping cone
$\mathsf{cone}(w)[-1]$ is $\tau_\kA$--periodic  and
isomorphic to a stalk complex $\kX[0]$, where $\kX$ is a locally projective $\kA$--module. Applying the functor $\TT$, we obtain a distinguished triangle
$$
\TT(\kX) \lar \mS_1[-1] \xrightarrow{\TT(w)} \mS_3[1] \lar \TT(\kX)[1].
$$
Next, we have:
$
\Hom_{D^b(\Lambda)}\bigl(\mS_1[-1], \mS_3[1]\bigr) \cong
\Ext^2_\Lambda(\mS_1, \mS_3) = \kk^2.
$
Note that $\mS_1$ has the projective resolution
$$
0 \lar \mP_3^{\oplus 2}
\xrightarrow{
\left(
\begin{smallmatrix}
b & 0 \\
0 & d
\end{smallmatrix}
\right)} \mP_2^{\oplus 2}
\xrightarrow{
\left(
\begin{smallmatrix}
a & c
\end{smallmatrix}
\right)}
\mP_1 \lar \mS_1 \lar 0.
$$
Hence, any element  $w \in \Ext^2_\Lambda(\mS_1, \mS_3)$ is given by a morphism
of $\Lambda$--modules
$(\lambda \mathbbm{1}, \mu \mathbbm{1}): \mP_3^{\oplus 2} \rightarrow \mP_3$, where
$(\lambda, \mu) \in \kk^2$.  Moreover, one can check
that the complex $\mathsf{cone}(w)$ is $\tau_\Lambda$--periodic if and only if
$(\lambda: \mu) \notin \{0, \infty\}$.  In that case, $\mathsf{cone}(w)$ is isomorphic
to the complex $\mU^\bullet_\nu$, where $\nu = - \frac{\mu}{\lambda}$.

\end{proof}

\section{Some generalizations and concluding remarks}

\subsection{Tilting on other degenerations of elliptic curves}

Theorem \ref{T:Main} can be generalized to the case of curves with more complicated singularities.
For example, let $\sE \subseteq \mathbb{P}^2$
be a tachnode plane cubic curve given by the equation $y (yz - x^2) = 0$.
Again, let $\pi: \widetilde\sE \rightarrow \sE$ be the normalization,
$\widetilde\kO = \pi_*(\kO_{\widetilde\sE})$ and $\kI = {\mathcal Ann}_\sE(\widetilde\kO/\kO)$ be the
conductor ideal.

 In this case,
consider again the coherent sheaf $\kF = \kI \oplus \kO$ and the sheaf of $\kO$--orders
$\kA = {\mathcal End}_\sE(\kF)$.
Let $\kA = \kP_1 \oplus \kP_2 \oplus \kF$ be the decomposition of $\kA$ into a direct
sum of indecomposable locally projective  modules. Define the torsion $\kA$--module $\kS$ via the short exact sequence
$$
0 \lar \kP_1 \oplus \kP_2 \lar  \kF \lar  \kS \lar 0.
$$
Then the complex  $\kH^\bullet = \kS[-1] \oplus \bigl(\kP_1 (-1)
\oplus \kP_1\bigr) \oplus \bigl(\kP_2(-1) \oplus \kP_2\bigr)$ is
rigid  and the endomorphism algebra $\Gamma_\sE = \End_{D^b(\kA)}(\kH^\bullet)$
is isomorphic to the path algebra
$$
\xymatrix
{
\bullet  &  \bullet  \ar@/^/[l]^{u_2}  \ar@/_/[l]_{v_2} &
 \bullet \ar[l]_{a_2} \ar[r]^{a_1}   \ar@(ul, ur)^\varphi  &
\bullet  \ar@/^/[r]^{u_1}  \ar@/_/[r]_{v_1} & \bullet
}
$$
subject to  the relations $\varphi^2 = 0, v_1 a_1 = 0, v_2 a_2 = 0, u_1 a_1 \varphi = 0$ and
$u_2 a_2 \varphi = 0$. Note  that in this case $\gldim(\Gamma_\sE) = \infty$. However, similar to the proof
of Theorem \ref{T:Main} one can show the
complex $\kH^\bullet$ is tilting in the sense of Theorem \ref{T:Keller}. Thus,
we have a fully faithful functor
$
\Perf(\sE) \rightarrow   D^b(\mod-\Gamma_\sE).
$

\subsection{Tilting on chains of projective lines}

Let $\sX$ be a chain of projective lines. In \cite{TorontoWorksh} it was shown that
$D\bigl(\Qcoh(\sX)\bigr)$ has a tilting vector bundle. In particular, we have
a triangle equivalence $D^*\bigl(\Coh(\sX)\bigr) \stackrel{\cong}\lar D^*(\rA_\sX-\mod)$,
where $* \in \{-,b \}$ and $\rA_\sX$ is the opposite algebra of the corresponding tilted algebra.
Composing this functor with the embedding obtained in Corollary \ref{C:embedofcategor},
get
an interesting fully faithful functor
$$
D^-(\rA_\sX-\mod) \lar D^-(\Lambda_\sX-\mod)
$$
which is worth to study in further details. Note that both algebras $\rA_\sX$ and $\Lambda_\sX$ are
gentle, hence derived-tame.  We hope that these geometric realizations of gentle algebras
will contribute to a better understanding of the representation theory of gentle algebras.

\begin{example}
Let $\sX = V(x_0 x_1) \subset \PP^2$ be a chain of two projective lines. Then $\rA_\sX$ is the path algebra
of the following quiver
$$
\xymatrix{ & \bullet & \\
\bullet \ar[ur]^a  \ar@/^/[rr]^c  & & \bullet \ar@/^/[ll]^d \ar[ul]_b
}
$$
subject to the relations $cd = 0 = dc$. The algebra $\Lambda_\sX$ is the path algebra of another
quiver
$$
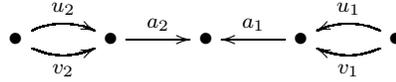
\xymatrix
{
\bullet  \ar@/^/[r]^{u_2}  \ar@/_/[r]_{v_2} &  \bullet  \ar[r]^{a_2} &
 \bullet   &
\bullet    \ar[l]_{a_1}& \bullet \ar@/^/[l]^{v_1}  \ar@/_/[l]_{u_1}
}
$$
subject to the relations $u_i a_i = 0$, $i = 1,2$.
\end{example}

\subsection{Non-commutative curves with nodal singularities}

Similar to constructions   of Section \ref{S:AuslanderSheaf} and Section \ref{S:Tilting},
one can study   a new class of derived-tame
non-commutative curves which generalizes weighted projective lines of Geigle and Lenzing
\cite{GeigleLenzing}.   Their characteristic property is that the completion of their
stalks are  generically matrix algebras over $\kk\llbracket t \rrbracket$, whereas
the \emph{singular}  stalks are either \emph{hereditary orders}  or \emph{nodal algebras}. The last
class of $\kk\llbracket t \rrbracket$--orders was
introduced in  \cite{Nodal}.

\begin{example}
Let $X = \mathbb{P}^1$,
 $Z = \bigl\{0, \infty\bigr\}$ and  $\kI = \kI_Z$ be the ideal sheaf of $Z$. Consider
the following sheaf of $\kO_{\mathbb{P}^1}$--orders:
$$
\kA =
\left(
\begin{array}{ccc}
\kO & \kI & \kI \\
\kO & \kO & \kI \\
\kO & \kI & \kO
\end{array}
\right).
$$
Then for $x \ne 0, \infty$ the algebra $\widehat{\kA_x} \cong
\Mat_3\bigl(\kk\llbracket t \rrbracket\bigr)$, whereas for $x \in \bigl\{0, \infty\bigr\}$
it is the completion of the
path algebra of the so-called \emph{Gelfand quiver}:
$$
G =
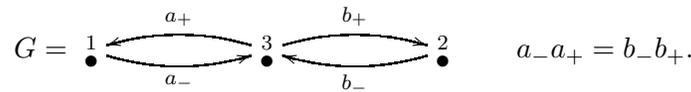
\xymatrix
{
\stackrel{1}\bullet \ar@/_/[rr]_{a_{-}}  & &  \stackrel{3}\bullet  \ar@/_/[ll]_{a_{+}}
 \ar@/^/[rr]^{b_{+}}
 & &
\ar@/^/[ll]^{b_{-}} \stackrel{2}\bullet}  \qquad a_{-} a_{+} = b_{-} b_{+}.
$$
Let $\mS_i$ be the simple $G$--module corresponding to the vertex $i  = 1,2$
and $\kS_i^{x}$ be the corresponding $\kA$--module for $x \in \left\{0, \infty \right\}$.
Let $e_1 \in H^0(\kA)$ be the primitive idempotent corresponding to the left upper corner unit element and
$$
\kP = \kA \cdot e_1 =
\left(
\begin{array}{c}
\kO \\
\kO \\
\kO
\end{array}
\right)
$$
be the corresponding locally projective $\kA$--module.
Similar to the proof of Proposition \ref{T:crucial}  and Theorem \ref{T:Main}
one can show that the complex
$\kH^\bullet = \bigl(\kS_1^0 \oplus \kS_2^0 \oplus \kS_1^\infty
\oplus \kS_2^\infty\bigr)[-1] \oplus
\kP(-1) \oplus \kP$  is tilting in $D^b\bigl(\Coh(\kA)\bigr)$. Its endomorphism algebra
$\Gamma  = \End_{D^b(\kA)}(\kH^\bullet)$ is isomorphic to the path algebra of
the quiver
$$
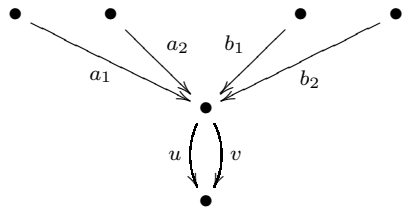
\xymatrix
{
\bullet \ar[rrd]_{a_1} & \bullet \ar[rd]^{a_2} & & \bullet \ar[ld]_{b_1} & \bullet \ar[lld]^{b_2} \\
                 &                 & \bullet \ar@/^/[d]^v  \ar@/_/[d]_u & & \\
                 &                 & \bullet  & &
}
$$
subject to the relations $u a_i = 0, \, v b_i = 0, \, i = 1,2$.
The algebra $\Gamma$ is derived equivalent to the path algebra of the quiver
$$
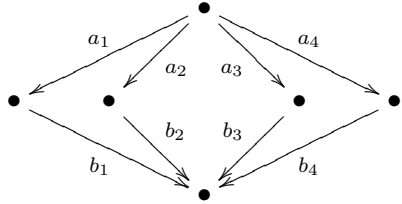
\xymatrix
{
        &           & \bullet \ar[lld]_{a_1}
\ar[ld]^{a_2}  \ar[rd]_{a_3}  \ar[rrd]^{a_4}
      &         &       \\
\bullet \ar[rrd]_{b_1}  & \bullet \ar[rd]^{b_2}
&        & \bullet \ar[ld]_{b_3}
& \bullet \ar[lld]^{b_4}\\
        &           &\bullet &         &       \\
}
$$
subject to the  relations $b_1 a_1 = b_2 a_2$ and  $b_3 a_3 = b_4 a_4$. This algebra
is a degeneration of the canonical tubular algebra $T_\lambda= T(2,2,2,2;\lambda)$,
$\lambda \in \kk \setminus\{0,1\}$,
introduced by Ringel in \cite{Ringel}.

The algebra $T(2,2,2,2;\lambda)$ has the following geometric interpretation. Consider
the elliptic curve $\sE \subseteq \PP^2$ given by the equation
$zy^2 = x(x-z)(x-\lambda z)$, where $\lambda \in \kk \setminus\{0,1\}$.
Then the group $\ZZ_2$ acts on $\sE$ by the rule $(x:y:z) \mapsto (x:-y:z)$.
By a result of Geigle and Lenzing \cite[Proposition 4.1 and Example 5.8]{GeigleLenzing},  see also
\cite[Corollary 1.4]{Orbifolds}, there is a derived equivalence
$$
D^b\bigl(\Coh^{\ZZ_2}(\sE)\bigr) \stackrel{\cong}\lar D^b(\mod-T_\lambda).
$$
Hence, the derived category $D^b\bigl(\Coh(\kA)\bigr)$ can be viewed as a \emph{degeneration}
of $D^b\bigl(\Coh^{\ZZ_2}(\sE)\bigr)$.
\end{example}

\subsection{Configuration schemes  of Lunts}

Our construction of non-commutative curves attached to a nodal rational projective curve,
is closely related with the category of coherent sheaves on a \emph{configuration scheme},
introduced by Lunts in \cite{Lunts}. Let $\sX$ be a union of projective lines intersecting transversally.
Lunts introduces a category $\Coh(\kX)$, constructs a fully faithful
exact functor
$\Perf(\sX) \rightarrow D^b\bigl(\Coh(\kX)\bigr)$ and shows that
$D^b\bigl(\Coh(\kX)\bigr)$
has a tilting \emph{sheaf}, whose endomorphism algebra is isomorphic to the algebra $\Gamma_\sX$ introduced
in Definition \ref{D:tiltedalg}.
Moreover, his approach can be generalized to higher dimensions, in particular to the case
of the singular quintic threefold $V(x_0 x_1 x_3 x_3 x_4) \subseteq \PP^4$ important in the string theory.

The  construction of Lunts seems  not to  generalize straightforwardly to the case of arbitrary
nodal rational projective curves.
Moreover, one can check that for a union of projective lines $\sX$ intersecting transversally,
the category $\Coh(\kX)$ is \emph{not}  equivalent to $\Coh(\kA_\sX)$.   In other words, $\Coh(\kX)$
is the heart of an interesting t-structure in the triangulated category
$D^b\bigl(\Coh(\kA_\sX)\bigr)$. An exact relation between the categories
$\Coh(\kX)$ and $\Coh(\kA_\sX)$  will be studied
in a separate paper.

\end{document}